\documentclass[12pt]{article}

\def\endofps{EndOfTheIncludedPostscriptMagicCookie}
\chardef\other=12
\newwrite\psdumphandle 
\outer\def\psdump#1{\par\medbreak
  \immediate\openout\psdumphandle=#1
  \copytoblankline}
\def\copytoblankline{\begingroup\setupcopy\copypsline}
\def\setupcopy{\def\do##1{\catcode`##1=\other}\dospecials
  \catcode`\\=\other \obeylines}
{\obeylines \gdef\copypsline#1
  {\def\next{#1}%
  \ifx\next\endofps\let\next=\endgroup %
  \else\immediate\write\psdumphandle{\next} \let\next=\copypsline\fi\next}}
\outer\def\closepsdump{
  \immediate\closeout\psdumphandle}

\psdump{altnewcover.eps}%!PS-Adobe-2.0 EPSF-2.0
/$F2psDict 200 dict def
$F2psDict begin
$F2psDict /mtrx matrix put
/col-1 {0 setgray} bind def
/col0 {0.000 0.000 0.000 srgb} bind def
/col1 {0.000 0.000 1.000 srgb} bind def
/col2 {0.000 1.000 0.000 srgb} bind def
/col3 {0.000 1.000 1.000 srgb} bind def
/col4 {1.000 0.000 0.000 srgb} bind def
/col5 {1.000 0.000 1.000 srgb} bind def
/col6 {1.000 1.000 0.000 srgb} bind def
/col7 {1.000 1.000 1.000 srgb} bind def
/col8 {0.000 0.000 0.560 srgb} bind def
/col9 {0.000 0.000 0.690 srgb} bind def
/col10 {0.000 0.000 0.820 srgb} bind def
/col11 {0.530 0.810 1.000 srgb} bind def
/col12 {0.000 0.560 0.000 srgb} bind def
/col13 {0.000 0.690 0.000 srgb} bind def
/col14 {0.000 0.820 0.000 srgb} bind def
/col15 {0.000 0.560 0.560 srgb} bind def
/col16 {0.000 0.690 0.690 srgb} bind def
/col17 {0.000 0.820 0.820 srgb} bind def
/col18 {0.560 0.000 0.000 srgb} bind def
/col19 {0.690 0.000 0.000 srgb} bind def
/col20 {0.820 0.000 0.000 srgb} bind def
/col21 {0.560 0.000 0.560 srgb} bind def
/col22 {0.690 0.000 0.690 srgb} bind def
/col23 {0.820 0.000 0.820 srgb} bind def
/col24 {0.500 0.190 0.000 srgb} bind def
/col25 {0.630 0.250 0.000 srgb} bind def
/col26 {0.750 0.380 0.000 srgb} bind def
/col27 {1.000 0.500 0.500 srgb} bind def
/col28 {1.000 0.630 0.630 srgb} bind def
/col29 {1.000 0.750 0.750 srgb} bind def
/col30 {1.000 0.880 0.880 srgb} bind def
/col31 {1.000 0.840 0.000 srgb} bind def

end
save
-5.0 497.0 translate
1 -1 scale

/cp {closepath} bind def
/ef {eofill} bind def
/gr {grestore} bind def
/gs {gsave} bind def
/sa {save} bind def
/rs {restore} bind def
/l {lineto} bind def
/m {moveto} bind def
/rm {rmoveto} bind def
/n {newpath} bind def
/s {stroke} bind def
/sh {show} bind def
/slc {setlinecap} bind def
/slj {setlinejoin} bind def
/slw {setlinewidth} bind def
/srgb {setrgbcolor} bind def
/rot {rotate} bind def
/sc {scale} bind def
/sd {setdash} bind def
/ff {findfont} bind def
/sf {setfont} bind def
/scf {scalefont} bind def
/sw {stringwidth} bind def
/tr {translate} bind def
/tnt {dup dup currentrgbcolor
  4 -2 roll dup 1 exch sub 3 -1 roll mul add
  4 -2 roll dup 1 exch sub 3 -1 roll mul add
  4 -2 roll dup 1 exch sub 3 -1 roll mul add srgb}
  bind def
/shd {dup dup currentrgbcolor 4 -2 roll mul 4 -2 roll mul
  4 -2 roll mul srgb} bind def
 /DrawEllipse {
/endangle exch def
/startangle exch def
/yrad exch def
/xrad exch def
/y exch def
/x exch def
/savematrix mtrx currentmatrix def
x y tr xrad yrad sc 0 0 1 startangle endangle arc
closepath
savematrix setmatrix
} def

/$F2psBegin {$F2psDict begin /$F2psEnteredState save def} def
/$F2psEnd {$F2psEnteredState restore end} def

$F2psBegin
10 setmiterlimit
n 0 792 m 0 0 l 612 0 l 612 792 l cp clip
 0.06000 0.06000 sc
1 slj
1 slc
15.000 slw
 [100.0] 0 sd
n 5803 5295 m 6361 5319 l gs col-1 s gr  [] 0 sd
n 6069 5191 m 5978 5422 l gs col-1 s gr
n 6171 5191 m 6080 5422 l gs col-1 s gr
 [100.0] 0 sd
n 1573 189 m 1984 249 l gs col-1 s gr  [] 0 sd
n 1768 317 m 1763 105 l gs col-1 s gr
 [100.0] 0 sd
n 1573 5285 m 1984 5345 l gs col-1 s gr  [] 0 sd
n 1768 5413 m 1763 5201 l gs col-1 s gr
n 6069 97 m 5978 328 l gs col-1 s gr
n 6171 97 m 6080 328 l gs col-1 s gr
 [100.0] 0 sd
n 5803 201 m 6361 225 l gs col-1 s gr  [] 0 sd
7.500 slw
n 8972 7119 48 48 0 360 DrawEllipse gs col7 0.00 shd ef gr gs col-1 s gr

15.000 slw
n 8825 6818 m 8975 6968 l 9125 6818 l gs col-1 s gr
 [100.0] 0 sd
n 5769 7120 m 6369 7147 l gs col-1 s gr  [] 0 sd
 [100.0] 0 sd
n 1792 7492 m 1672 7957 l gs col-1 s gr  [] 0 sd
30.000 slw
gs n 5760 6313 m
6624.2 6308.2 6997.4 6312.9 7253 6332 curveto
7605.5 6358.3 8436.4 6440.1 8783 6542 curveto
9049.1 6620.2 9808.7 6690.9 9858 7042 curveto
9899.5 7337.2 9305.4 7572.8 9105 7695 curveto
8910.1 7813.8 8418.6 7991.1 8205 8055 curveto
7989.1 8119.5 7478.7 8249.4 7245 8250 curveto
6949.6 8250.8 6315.3 8034.0 6030 8010 curveto
5774.2 7988.5 5175.8 8061.3 4920 8040 curveto
4633.2 8016.1 3988.3 7815.6 3698 7801 curveto
3489.4 7790.5 3009.9 7864.6 2807 7884 curveto
2670.6 7897.0 2471.1 7916.5 2009 7962 curveto
 gs col-1 s gr
 gr

gs n 7327 6851 m
7432.9 6845.8 7684.1 6837.1 7791 6851 curveto
7893.2 6864.2 8124.9 6898.5 8220 6964 curveto
8269.3 6998.0 8362.3 7067.5 8354 7151 curveto
8342.9 7262.7 8193.8 7324.0 8129 7361 curveto
8042.5 7410.4 7828.6 7460.8 7735 7477 curveto
7647.5 7492.2 7440.0 7500.9 7352 7499 curveto
7227.0 7496.3 6935.1 7488.6 6810 7444 curveto
6688.1 7400.5 6362.4 7304.8 6367 7151 curveto
6371.2 7007.8 6681.8 6939.4 6795 6903 curveto
6917.3 6863.6 7205.0 6856.9 7327 6851 curveto
 cp  gs col-1 s gr
 gr

gs n 5010 6784 m
4912.1 6790.0 4700.4 6818.8 4612 6889 curveto
4532.1 6952.4 4402.1 7083.4 4438 7218 curveto
4477.3 7365.6 4681.9 7410.7 4783 7436 curveto
4884.6 7461.4 5097.9 7431.1 5192 7406 curveto
5338.4 7366.9 5712.2 7347.0 5767 7129 curveto
5793.2 7024.6 5654.7 6926.0 5587 6882 curveto
5531.1 6845.7 5396.8 6838.9 5340 6829 curveto
5265.1 6815.9 5089.3 6779.1 5010 6784 curveto
 cp  gs col-1 s gr
 gr

gs n 2009 7962 m
1836.5 7959.6 1762.0 7956.6 1711 7950 curveto
1517.0 7925.1 1066.0 7836.5 880 7768 curveto
716.0 7707.6 347.3 7567.9 218 7415 curveto
160.4 7346.9 106.3 7187.5 112 7098 curveto
116.5 7027.2 186.6 6909.8 227 6860 curveto
293.4 6778.2 474.4 6637.9 571 6595 curveto
912.1 6443.5 1717.7 6322.5 2080 6321 curveto
2312.6 6320.1 2821.8 6478.2 3050 6489 curveto
3210.2 6496.6 3582.3 6454.4 3739 6436 curveto
3910.9 6415.8 4314.6 6331.1 4488 6317 curveto
4706.1 6299.3 5024.1 6297.8 5760 6311 curveto
 gs col-1 s gr
 gr

gs n 1954 7495 m
1792.2 7499.9 1722.2 7499.9 1674 7495 curveto
1523.1 7479.6 1174.4 7422.3 1030 7362 curveto
932.1 7321.1 654.9 7244.3 659 7115 curveto
664.3 6946.4 1027.3 6864.9 1162 6824 curveto
1410.1 6748.6 1978.0 6723.8 2230 6754 curveto
2341.0 6767.3 2582.0 6827.6 2680 6895 curveto
2743.6 6938.7 2880.5 7019.4 2870 7135 curveto
2857.3 7275.2 2674.0 7340.1 2591 7380 curveto
2481.3 7432.8 2320.3 7461.5 1947 7495 curveto
 gs col-1 s gr
 gr

7.500 slw
n 8972 3997 48 48 0 360 DrawEllipse gs col7 0.00 shd ef gr gs col-1 s gr

n 8972 1425 48 48 0 360 DrawEllipse gs col7 0.00 shd ef gr gs col-1 s gr

15.000 slw
n 8974 1421 m 8975 6968 l gs col-1 s gr
30.000 slw
gs n 1460 3028 m
1173.5 3009.6 1048.7 3009.6 961 3028 curveto
822.0 3057.1 531.9 3173.0 418 3266 curveto
326.1 3341.1 153.9 3512.0 142 3657 curveto
132.3 3775.9 233.5 3939.8 313 4010 curveto
408.9 4094.7 637.0 4165.5 761 4162 curveto
887.6 4158.4 1113.5 4049.6 1209 3981 curveto
1315.8 3904.3 1488.2 3671.6 1561 3571 curveto
1632.9 3471.7 1780.2 3222.1 1837 3114 curveto
1864.7 3061.3 1900.5 2982.3 1980 2798 curveto
 gs col-1 s gr
 gr

gs n 6857 4237 m
7015.6 3980.5 7098.9 3877.5 7190 3825 curveto
7318.1 3751.2 7605.8 3726.0 7742 3739 curveto
7904.4 3754.5 8308.9 3757.4 8369 3964 curveto
8401.1 4074.3 8224.8 4193.6 8151 4245 curveto
8059.8 4308.5 7832.6 4359.2 7732 4377 curveto
7633.9 4394.3 7403.1 4404.0 7304 4396 curveto
7223.6 4389.5 7041.2 4375.2 6964 4325 curveto
6720.3 4166.5 6403.9 3707.4 6324 3435 curveto
6287.9 3312.0 6296.7 3147.2 6359 2776 curveto
 gs col-1 s gr
 gr

gs n 5833 4976 m
5825.6 5065.2 5822.3 5103.7 5820 5130 curveto
5817.6 5157.5 5813.8 5197.8 5805 5291 curveto
 gs col-1 s gr
 gr

gs n 6428 5080 m
6402.5 5147.9 6392.0 5177.4 6386 5198 curveto
6380.5 5216.9 6374.0 5244.9 6360 5310 curveto
 gs col-1 s gr
 gr

gs n 1863 3061 m
2619.9 3157.5 2947.4 3194.0 3173 3207 curveto
3473.1 3224.2 4188.3 3216.4 4488 3212 curveto
4702.9 3208.8 5017.2 3198.3 5745 3170 curveto
 gs col-1 s gr
 gr

gs n 6306 3170 m
6854.6 3196.8 7091.4 3211.1 7253 3227 curveto
7603.7 3261.6 8436.4 3335.1 8783 3437 curveto
9049.1 3515.2 9820.3 3600.1 9858 3937 curveto
9895.1 4269.2 9172.7 4519.6 8933 4652 curveto
8757.2 4749.1 8319.0 4897.0 8128 4947 curveto
7842.5 5021.7 7165.2 5155.1 6863 5152 curveto
6647.2 5149.8 6170.0 5015.5 5963 4987 curveto
5644.9 4943.2 4882.6 4880.8 4563 4857 curveto
4298.9 4837.4 3674.0 4765.2 3403 4792 curveto
3091.6 4822.8 2397.6 4914.7 2110 5126 curveto
2069.3 5155.9 2036.1 5209.4 1977 5340 curveto
 gs col-1 s gr
 gr

gs n 5805 2759 m
5745.0 3225.7 5708.5 3426.4 5659 3562 curveto
5592.3 3744.7 5451.0 4148.9 5249 4278 curveto
5126.6 4356.2 4876.3 4353.2 4751 4321 curveto
4656.3 4296.6 4485.5 4217.7 4437 4111 curveto
4398.8 4027.0 4419.5 3887.4 4475 3816 curveto
4581.0 3679.7 4844.3 3597.1 5003 3607 curveto
5127.1 3614.8 5256.9 3686.0 5522 3892 curveto
 gs col-1 s gr
 gr

gs n 1312 3355 m
1110.8 3343.0 1023.1 3344.5 961 3361 curveto
881.7 3382.1 716.9 3430.8 666 3523 curveto
642.8 3565.0 632.9 3636.3 675 3676 curveto
766.8 3762.5 926.3 3708.7 1009 3676 curveto
1090.8 3643.6 1202.7 3511.6 1250 3450 curveto
1299.2 3385.9 1378.9 3218.2 1412 3147 curveto
1443.8 3078.5 1484.6 2976.3 1575 2738 curveto
 gs col-1 s gr
 gr

gs n 1679 3392 m
2092.3 3411.8 2271.1 3433.3 2394 3478 curveto
2497.0 3515.4 2701.8 3631.5 2773 3725 curveto
2823.4 3791.2 2869.0 3949.4 2870 4030 curveto
2870.9 4107.3 2830.1 4260.3 2784 4325 curveto
2721.0 4413.5 2525.8 4520.7 2441 4568 curveto
2364.3 4610.7 2168.9 4677.3 2091 4716 curveto
2019.8 4751.3 1850.4 4827.7 1790 4889 curveto
1732.4 4947.5 1678.9 5044.0 1576 5275 curveto
 gs col-1 s gr
 gr

gs n 1460 485 m
1173.5 466.6 1048.7 466.6 961 485 curveto
822.0 514.1 531.9 630.0 418 723 curveto
326.1 798.0 153.9 969.0 142 1114 curveto
132.3 1232.9 233.5 1396.8 313 1467 curveto
408.9 1551.7 637.0 1622.5 761 1619 curveto
887.6 1615.4 1113.5 1506.6 1209 1438 curveto
1315.8 1361.3 1488.2 1128.6 1561 1028 curveto
1632.9 928.6 1780.2 679.1 1837 571 curveto
1864.7 518.3 1900.5 439.3 1980 255 curveto
 gs col-1 s gr
 gr

gs n 6857 1694 m
7015.6 1437.5 7098.9 1334.5 7190 1282 curveto
7318.1 1208.2 7605.8 1183.0 7742 1196 curveto
7904.4 1211.5 8308.9 1214.4 8369 1421 curveto
8401.1 1531.3 8224.8 1650.6 8151 1702 curveto
8059.8 1765.5 7832.6 1816.2 7732 1834 curveto
7633.9 1851.3 7403.1 1861.0 7304 1853 curveto
7223.6 1846.5 7041.2 1832.2 6964 1782 curveto
6720.3 1623.5 6403.9 1164.4 6324 892 curveto
6287.9 769.0 6296.7 604.2 6359 233 curveto
 gs col-1 s gr
 gr

gs n 5833 2433 m
5825.6 2522.2 5822.3 2560.7 5820 2587 curveto
5817.6 2614.5 5813.8 2654.8 5805 2748 curveto
 gs col-1 s gr
 gr

gs n 6428 2537 m
6402.5 2604.9 6392.0 2634.4 6386 2655 curveto
6380.5 2673.9 6374.0 2701.9 6360 2767 curveto
 gs col-1 s gr
 gr

gs n 1863 518 m
2619.9 614.5 2947.4 651.0 3173 664 curveto
3473.1 681.2 4188.3 673.4 4488 669 curveto
4702.9 665.8 5017.2 655.3 5745 627 curveto
 gs col-1 s gr
 gr

gs n 6306 627 m
6854.6 653.8 7091.4 668.1 7253 684 curveto
7603.7 718.6 8436.4 792.1 8783 894 curveto
9049.1 972.2 9820.3 1057.1 9858 1394 curveto
9895.1 1726.2 9172.7 1976.6 8933 2109 curveto
8757.2 2206.1 8319.0 2354.0 8128 2404 curveto
7842.5 2478.7 7165.2 2612.1 6863 2609 curveto
6647.2 2606.8 6170.0 2472.5 5963 2444 curveto
5644.9 2400.2 4882.6 2337.8 4563 2314 curveto
4298.9 2294.3 3674.0 2222.2 3403 2249 curveto
3091.6 2279.8 2397.6 2371.7 2110 2583 curveto
2069.3 2612.9 2036.1 2666.4 1977 2797 curveto
 gs col-1 s gr
 gr

gs n 5805 216 m
5745.0 682.7 5708.5 883.4 5659 1019 curveto
5592.3 1201.7 5451.0 1606.0 5249 1735 curveto
5126.6 1813.2 4876.3 1810.2 4751 1778 curveto
4656.3 1753.7 4485.5 1674.7 4437 1568 curveto
4398.8 1484.0 4419.5 1344.4 4475 1273 curveto
4581.0 1136.7 4844.3 1054.0 5003 1064 curveto
5127.1 1071.8 5256.9 1143.0 5522 1349 curveto
 gs col-1 s gr
 gr

gs n 1312 812 m
1110.8 800.0 1023.1 801.5 961 818 curveto
881.7 839.1 716.9 887.8 666 980 curveto
642.8 1022.0 632.9 1093.3 675 1133 curveto
766.8 1219.5 926.3 1165.7 1009 1133 curveto
1090.8 1100.6 1202.7 968.6 1250 907 curveto
1299.2 842.9 1378.9 675.2 1412 604 curveto
1443.8 535.5 1484.6 433.3 1575 195 curveto
 gs col-1 s gr
 gr

gs n 1679 849 m
2092.3 868.8 2271.1 890.3 2394 935 curveto
2497.0 972.4 2701.8 1088.5 2773 1182 curveto
2823.4 1248.2 2869.0 1406.5 2870 1487 curveto
2870.9 1564.3 2830.1 1717.3 2784 1782 curveto
2721.0 1870.5 2525.8 1977.7 2441 2025 curveto
2364.3 2067.7 2168.9 2134.3 2091 2173 curveto
2019.8 2208.3 1850.4 2284.7 1790 2346 curveto
1732.4 2404.5 1678.9 2501.0 1576 2732 curveto
 gs col-1 s gr
 gr

$F2psEnd
rs
EndOfTheIncludedPostscriptMagicCookie

\closepsdump

\psdump{klines.eps}%!PS-Adobe-2.0 EPSF-2.0
/$F2psDict 200 dict def
$F2psDict begin
$F2psDict /mtrx matrix put
/col-1 {0 setgray} bind def
/col0 {0.000 0.000 0.000 srgb} bind def
/col1 {0.000 0.000 1.000 srgb} bind def
/col2 {0.000 1.000 0.000 srgb} bind def
/col3 {0.000 1.000 1.000 srgb} bind def
/col4 {1.000 0.000 0.000 srgb} bind def
/col5 {1.000 0.000 1.000 srgb} bind def
/col6 {1.000 1.000 0.000 srgb} bind def
/col7 {1.000 1.000 1.000 srgb} bind def
/col8 {0.000 0.000 0.560 srgb} bind def
/col9 {0.000 0.000 0.690 srgb} bind def
/col10 {0.000 0.000 0.820 srgb} bind def
/col11 {0.530 0.810 1.000 srgb} bind def
/col12 {0.000 0.560 0.000 srgb} bind def
/col13 {0.000 0.690 0.000 srgb} bind def
/col14 {0.000 0.820 0.000 srgb} bind def
/col15 {0.000 0.560 0.560 srgb} bind def
/col16 {0.000 0.690 0.690 srgb} bind def
/col17 {0.000 0.820 0.820 srgb} bind def
/col18 {0.560 0.000 0.000 srgb} bind def
/col19 {0.690 0.000 0.000 srgb} bind def
/col20 {0.820 0.000 0.000 srgb} bind def
/col21 {0.560 0.000 0.560 srgb} bind def
/col22 {0.690 0.000 0.690 srgb} bind def
/col23 {0.820 0.000 0.820 srgb} bind def
/col24 {0.500 0.190 0.000 srgb} bind def
/col25 {0.630 0.250 0.000 srgb} bind def
/col26 {0.750 0.380 0.000 srgb} bind def
/col27 {1.000 0.500 0.500 srgb} bind def
/col28 {1.000 0.630 0.630 srgb} bind def
/col29 {1.000 0.750 0.750 srgb} bind def
/col30 {1.000 0.880 0.880 srgb} bind def
/col31 {1.000 0.840 0.000 srgb} bind def

end
save
-34.0 273.0 translate
1 -1 scale

/cp {closepath} bind def
/ef {eofill} bind def
/gr {grestore} bind def
/gs {gsave} bind def
/sa {save} bind def
/rs {restore} bind def
/l {lineto} bind def
/m {moveto} bind def
/rm {rmoveto} bind def
/n {newpath} bind def
/s {stroke} bind def
/sh {show} bind def
/slc {setlinecap} bind def
/slj {setlinejoin} bind def
/slw {setlinewidth} bind def
/srgb {setrgbcolor} bind def
/rot {rotate} bind def
/sc {scale} bind def
/sd {setdash} bind def
/ff {findfont} bind def
/sf {setfont} bind def
/scf {scalefont} bind def
/sw {stringwidth} bind def
/tr {translate} bind def
/tnt {dup dup currentrgbcolor
  4 -2 roll dup 1 exch sub 3 -1 roll mul add
  4 -2 roll dup 1 exch sub 3 -1 roll mul add
  4 -2 roll dup 1 exch sub 3 -1 roll mul add srgb}
  bind def
/shd {dup dup currentrgbcolor 4 -2 roll mul 4 -2 roll mul
  4 -2 roll mul srgb} bind def
 /DrawEllipse {
/endangle exch def
/startangle exch def
/yrad exch def
/xrad exch def
/y exch def
/x exch def
/savematrix mtrx currentmatrix def
x y tr xrad yrad sc 0 0 1 startangle endangle arc
closepath
savematrix setmatrix
} def

/$F2psBegin {$F2psDict begin /$F2psEnteredState save def} def
/$F2psEnd {$F2psEnteredState restore end} def

$F2psBegin
10 setmiterlimit
n 0 792 m 0 0 l 612 0 l 612 792 l cp clip
 0.06000 0.06000 sc
30.000 slw
n 2295 3270 229 229 0 360 DrawEllipse gs col-1 s gr

n 3900 3165 342 342 0 360 DrawEllipse gs col-1 s gr

7.500 slw
n 3847 3502 54 54 0 360 DrawEllipse gs col7 0.00 shd ef gr gs col-1 s gr

30.000 slw
n 3000 1800 361 361 0 360 DrawEllipse gs col-1 s gr

7.500 slw
n 3075 2152 54 54 0 360 DrawEllipse gs col7 0.00 shd ef gr gs col-1 s gr

n 2275 3502 54 54 0 360 DrawEllipse gs col7 0.00 shd ef gr gs col-1 s gr

n 2595 4477 54 54 0 360 DrawEllipse gs col7 0.00 shd ef gr gs col-1 s gr

30.000 slw
gs n 3375 675 m
3785.0 752.4 4580.8 1072.5 4875 1425 curveto
5203.9 1819.1 5598.7 2643.2 5400 3225 curveto
5183.8 3858.1 4291.3 4301.4 3750 4425 curveto
3082.5 4577.4 1854.3 4499.9 1275 3975 curveto
772.1 3519.3 508.6 2517.3 675 1875 curveto
788.1 1438.3 1374.5 989.2 1725 825 curveto
2116.1 641.8 2969.9 598.6 3375 675 curveto
 cp  gs col-1 s gr
 gr

1 slc
15.000 slw
gs n 3075 2145 m
3126.2 2392.3 3141.2 2501.0 3135 2580 curveto
3101.4 3010.0 2705.4 3895.2 2617 4297 curveto
2610.7 4325.9 2604.7 4368.4 2593 4467 curveto
 gs col-1 s gr
 gr

gs n 3855 3510 m
3805.7 3629.0 3781.9 3679.5 3760 3712 curveto
3689.6 3816.5 3497.9 4031.2 3395 4107 curveto
3275.2 4195.3 2961.7 4321.3 2830 4377 curveto
2790.0 4393.9 2730.5 4415.9 2592 4465 curveto
 gs col-1 s gr
 gr

gs n 2275 3497 m
2258.8 3709.9 2256.3 3802.4 2265 3867 curveto
2279.6 3975.3 2334.0 4214.8 2400 4312 curveto
2426.7 4351.3 2475.7 4390.3 2596 4468 curveto
 gs col-1 s gr
 gr

$F2psEnd
rs
EndOfTheIncludedPostscriptMagicCookie

\closepsdump

\psdump{contrad.eps}%!PS-Adobe-2.0 EPSF-2.0
/$F2psDict 200 dict def
$F2psDict begin
$F2psDict /mtrx matrix put
/col-1 {0 setgray} bind def
/col0 {0.000 0.000 0.000 srgb} bind def
/col1 {0.000 0.000 1.000 srgb} bind def
/col2 {0.000 1.000 0.000 srgb} bind def
/col3 {0.000 1.000 1.000 srgb} bind def
/col4 {1.000 0.000 0.000 srgb} bind def
/col5 {1.000 0.000 1.000 srgb} bind def
/col6 {1.000 1.000 0.000 srgb} bind def
/col7 {1.000 1.000 1.000 srgb} bind def
/col8 {0.000 0.000 0.560 srgb} bind def
/col9 {0.000 0.000 0.690 srgb} bind def
/col10 {0.000 0.000 0.820 srgb} bind def
/col11 {0.530 0.810 1.000 srgb} bind def
/col12 {0.000 0.560 0.000 srgb} bind def
/col13 {0.000 0.690 0.000 srgb} bind def
/col14 {0.000 0.820 0.000 srgb} bind def
/col15 {0.000 0.560 0.560 srgb} bind def
/col16 {0.000 0.690 0.690 srgb} bind def
/col17 {0.000 0.820 0.820 srgb} bind def
/col18 {0.560 0.000 0.000 srgb} bind def
/col19 {0.690 0.000 0.000 srgb} bind def
/col20 {0.820 0.000 0.000 srgb} bind def
/col21 {0.560 0.000 0.560 srgb} bind def
/col22 {0.690 0.000 0.690 srgb} bind def
/col23 {0.820 0.000 0.820 srgb} bind def
/col24 {0.500 0.190 0.000 srgb} bind def
/col25 {0.630 0.250 0.000 srgb} bind def
/col26 {0.750 0.380 0.000 srgb} bind def
/col27 {1.000 0.500 0.500 srgb} bind def
/col28 {1.000 0.630 0.630 srgb} bind def
/col29 {1.000 0.750 0.750 srgb} bind def
/col30 {1.000 0.880 0.880 srgb} bind def
/col31 {1.000 0.840 0.000 srgb} bind def

end
save
-3.0 225.0 translate
1 -1 scale

/cp {closepath} bind def
/ef {eofill} bind def
/gr {grestore} bind def
/gs {gsave} bind def
/sa {save} bind def
/rs {restore} bind def
/l {lineto} bind def
/m {moveto} bind def
/rm {rmoveto} bind def
/n {newpath} bind def
/s {stroke} bind def
/sh {show} bind def
/slc {setlinecap} bind def
/slj {setlinejoin} bind def
/slw {setlinewidth} bind def
/srgb {setrgbcolor} bind def
/rot {rotate} bind def
/sc {scale} bind def
/sd {setdash} bind def
/ff {findfont} bind def
/sf {setfont} bind def
/scf {scalefont} bind def
/sw {stringwidth} bind def
/tr {translate} bind def
/tnt {dup dup currentrgbcolor
  4 -2 roll dup 1 exch sub 3 -1 roll mul add
  4 -2 roll dup 1 exch sub 3 -1 roll mul add
  4 -2 roll dup 1 exch sub 3 -1 roll mul add srgb}
  bind def
/shd {dup dup currentrgbcolor 4 -2 roll mul 4 -2 roll mul
  4 -2 roll mul srgb} bind def
 /DrawEllipse {
/endangle exch def
/startangle exch def
/yrad exch def
/xrad exch def
/y exch def
/x exch def
/savematrix mtrx currentmatrix def
x y tr xrad yrad sc 0 0 1 startangle endangle arc
closepath
savematrix setmatrix
} def

/$F2psBegin {$F2psDict begin /$F2psEnteredState save def} def
/$F2psEnd {$F2psEnteredState restore end} def

$F2psBegin
10 setmiterlimit
n 0 792 m 0 0 l 612 0 l 612 792 l cp clip
 0.06000 0.06000 sc
30.000 slw
n 2280 1860 511 511 0 360 DrawEllipse gs col-1 s gr

7.500 slw
n 1777 1732 40 40 0 360 DrawEllipse gs col7 0.00 shd ef gr gs col-1 s gr

n 907 375 40 40 0 360 DrawEllipse gs col7 0.00 shd ef gr gs col-1 s gr

n 1190 1015 40 40 0 360 DrawEllipse gs col7 0.00 shd ef gr gs col-1 s gr

30.000 slw
gs n 2712 87 m
3099.6 160.2 3851.9 462.7 4130 796 curveto
4440.9 1168.3 4814.1 1947.1 4626 2498 curveto
4421.9 3095.3 3578.2 3514.2 3067 3631 curveto
2435.6 3775.1 1274.6 3701.8 727 3206 curveto
251.6 2775.1 2.5 1828.1 160 1221 curveto
266.8 808.4 821.0 383.9 1152 229 curveto
1522.1 55.7 2329.2 14.8 2712 87 curveto
 cp  gs col-1 s gr
 gr

1 slc
15.000 slw
gs  [100.0] 0 sd
n 1185 1005 m
1237.4 1318.0 1280.4 1451.0 1357 1537 curveto
1421.4 1609.3 1528.4 1658.0 1785 1732 curveto
 gs col-1 s gr
 gr
 [] 0 sd
gs n 925 385 m
1046.7 551.5 1095.4 626.5 1120 685 curveto
1142.2 737.8 1159.7 816.5 1190 1000 curveto
 gs col-1 s gr
 gr

gs n 920 385 m
991.6 484.9 1015.3 533.6 1015 580 curveto
1012.8 899.8 536.5 1284.0 510 1582 curveto
482.0 1895.9 645.5 2467.5 847 2715 curveto
1054.7 2970.1 1619.7 3231.4 1920 3292 curveto
2289.4 3366.6 3055.3 3337.0 3412 3150 curveto
3701.3 2998.3 4120.2 2566.8 4195 2220 curveto
4267.6 1883.6 4061.5 1334.8 3875 1085 curveto
3674.9 817.0 3126.5 481.1 2805 410 curveto
2525.7 348.3 1966.9 437.9 1715 550 curveto
1594.9 603.4 1462.4 717.2 1185 1005 curveto
 gs col-1 s gr
 gr

gs  [100.0] 0 sd
n 1190 1005 m
1010.7 1313.6 949.4 1457.9 945 1582 curveto
937.1 1801.9 1086.4 2200.4 1215 2370 curveto
1345.7 2542.3 1706.1 2787.9 1912 2850 curveto
2121.9 2913.3 2565.6 2912.4 2775 2835 curveto
2978.0 2760.0 3318.9 2511.6 3420 2302 curveto
3516.9 2101.1 3513.7 1701.3 3435 1500 curveto
3354.2 1293.3 3067.7 994.8 2865 905 curveto
2638.6 804.7 2189.3 746.6 1945 880 curveto
1743.4 990.1 1509.7 1262.0 1552 1530 curveto
1564.9 1611.6 1623.1 1662.1 1785 1732 curveto
 gs col-1 s gr
 gr
 [] 0 sd
$F2psEnd
rs
EndOfTheIncludedPostscriptMagicCookie

\closepsdump

\psdump{movex.eps}%!PS-Adobe-2.0 EPSF-2.0
/$F2psDict 200 dict def
$F2psDict begin
$F2psDict /mtrx matrix put
/col-1 {0 setgray} bind def
/col0 {0.000 0.000 0.000 srgb} bind def
/col1 {0.000 0.000 1.000 srgb} bind def
/col2 {0.000 1.000 0.000 srgb} bind def
/col3 {0.000 1.000 1.000 srgb} bind def
/col4 {1.000 0.000 0.000 srgb} bind def
/col5 {1.000 0.000 1.000 srgb} bind def
/col6 {1.000 1.000 0.000 srgb} bind def
/col7 {1.000 1.000 1.000 srgb} bind def
/col8 {0.000 0.000 0.560 srgb} bind def
/col9 {0.000 0.000 0.690 srgb} bind def
/col10 {0.000 0.000 0.820 srgb} bind def
/col11 {0.530 0.810 1.000 srgb} bind def
/col12 {0.000 0.560 0.000 srgb} bind def
/col13 {0.000 0.690 0.000 srgb} bind def
/col14 {0.000 0.820 0.000 srgb} bind def
/col15 {0.000 0.560 0.560 srgb} bind def
/col16 {0.000 0.690 0.690 srgb} bind def
/col17 {0.000 0.820 0.820 srgb} bind def
/col18 {0.560 0.000 0.000 srgb} bind def
/col19 {0.690 0.000 0.000 srgb} bind def
/col20 {0.820 0.000 0.000 srgb} bind def
/col21 {0.560 0.000 0.560 srgb} bind def
/col22 {0.690 0.000 0.690 srgb} bind def
/col23 {0.820 0.000 0.820 srgb} bind def
/col24 {0.500 0.190 0.000 srgb} bind def
/col25 {0.630 0.250 0.000 srgb} bind def
/col26 {0.750 0.380 0.000 srgb} bind def
/col27 {1.000 0.500 0.500 srgb} bind def
/col28 {1.000 0.630 0.630 srgb} bind def
/col29 {1.000 0.750 0.750 srgb} bind def
/col30 {1.000 0.880 0.880 srgb} bind def
/col31 {1.000 0.840 0.000 srgb} bind def

end
save
-34.0 273.0 translate
1 -1 scale

/cp {closepath} bind def
/ef {eofill} bind def
/gr {grestore} bind def
/gs {gsave} bind def
/sa {save} bind def
/rs {restore} bind def
/l {lineto} bind def
/m {moveto} bind def
/rm {rmoveto} bind def
/n {newpath} bind def
/s {stroke} bind def
/sh {show} bind def
/slc {setlinecap} bind def
/slj {setlinejoin} bind def
/slw {setlinewidth} bind def
/srgb {setrgbcolor} bind def
/rot {rotate} bind def
/sc {scale} bind def
/sd {setdash} bind def
/ff {findfont} bind def
/sf {setfont} bind def
/scf {scalefont} bind def
/sw {stringwidth} bind def
/tr {translate} bind def
/tnt {dup dup currentrgbcolor
  4 -2 roll dup 1 exch sub 3 -1 roll mul add
  4 -2 roll dup 1 exch sub 3 -1 roll mul add
  4 -2 roll dup 1 exch sub 3 -1 roll mul add srgb}
  bind def
/shd {dup dup currentrgbcolor 4 -2 roll mul 4 -2 roll mul
  4 -2 roll mul srgb} bind def
 /DrawEllipse {
/endangle exch def
/startangle exch def
/yrad exch def
/xrad exch def
/y exch def
/x exch def
/savematrix mtrx currentmatrix def
x y tr xrad yrad sc 0 0 1 startangle endangle arc
closepath
savematrix setmatrix
} def

/$F2psBegin {$F2psDict begin /$F2psEnteredState save def} def
/$F2psEnd {$F2psEnteredState restore end} def

$F2psBegin
10 setmiterlimit
n 0 792 m 0 0 l 612 0 l 612 792 l cp clip
 0.06000 0.06000 sc
30.000 slw
n 2295 3270 229 229 0 360 DrawEllipse gs col-1 s gr

n 3900 3165 342 342 0 360 DrawEllipse gs col-1 s gr

n 3000 1800 361 361 0 360 DrawEllipse gs col-1 s gr

7.500 slw
n 3075 2152 54 54 0 360 DrawEllipse gs col7 0.00 shd ef gr gs col-1 s gr

n 2595 4477 54 54 0 360 DrawEllipse gs col7 0.00 shd ef gr gs col-1 s gr

n 2790 2083 54 54 0 360 DrawEllipse gs col7 0.00 shd ef gr gs col-1 s gr

n 1675 4228 54 54 0 360 DrawEllipse gs col7 0.00 shd ef gr gs col-1 s gr

30.000 slw
gs n 3375 675 m
3785.0 752.4 4580.8 1072.5 4875 1425 curveto
5203.9 1819.1 5598.7 2643.2 5400 3225 curveto
5183.8 3858.1 4291.3 4301.4 3750 4425 curveto
3082.5 4577.4 1854.3 4499.9 1275 3975 curveto
772.1 3519.3 508.6 2517.3 675 1875 curveto
788.1 1438.3 1374.5 989.2 1725 825 curveto
2116.1 641.8 2969.9 598.6 3375 675 curveto
 cp  gs col-1 s gr
 gr

1 slc
15.000 slw
gs n 2785 2083 m
2713.2 2192.6 2685.7 2242.6 2675 2283 curveto
2600.6 2565.3 2839.1 3169.3 2660 3478 curveto
2507.8 3740.3 1958.5 3807.2 1785 3993 curveto
1752.0 4028.3 1723.3 4085.8 1670 4223 curveto
 gs col-1 s gr
 gr

gs n 3075 2145 m
3126.2 2392.3 3141.2 2501.0 3135 2580 curveto
3101.4 3010.0 2705.4 3895.2 2617 4297 curveto
2610.7 4325.9 2604.7 4368.4 2593 4467 curveto
 gs col-1 s gr
 gr

$F2psEnd
rs
EndOfTheIncludedPostscriptMagicCookie

\closepsdump

\psdump{cutting.eps}%!PS-Adobe-2.0 EPSF-2.0
/$F2psDict 200 dict def
$F2psDict begin
$F2psDict /mtrx matrix put
/col-1 {0 setgray} bind def
/col0 {0.000 0.000 0.000 srgb} bind def
/col1 {0.000 0.000 1.000 srgb} bind def
/col2 {0.000 1.000 0.000 srgb} bind def
/col3 {0.000 1.000 1.000 srgb} bind def
/col4 {1.000 0.000 0.000 srgb} bind def
/col5 {1.000 0.000 1.000 srgb} bind def
/col6 {1.000 1.000 0.000 srgb} bind def
/col7 {1.000 1.000 1.000 srgb} bind def
/col8 {0.000 0.000 0.560 srgb} bind def
/col9 {0.000 0.000 0.690 srgb} bind def
/col10 {0.000 0.000 0.820 srgb} bind def
/col11 {0.530 0.810 1.000 srgb} bind def
/col12 {0.000 0.560 0.000 srgb} bind def
/col13 {0.000 0.690 0.000 srgb} bind def
/col14 {0.000 0.820 0.000 srgb} bind def
/col15 {0.000 0.560 0.560 srgb} bind def
/col16 {0.000 0.690 0.690 srgb} bind def
/col17 {0.000 0.820 0.820 srgb} bind def
/col18 {0.560 0.000 0.000 srgb} bind def
/col19 {0.690 0.000 0.000 srgb} bind def
/col20 {0.820 0.000 0.000 srgb} bind def
/col21 {0.560 0.000 0.560 srgb} bind def
/col22 {0.690 0.000 0.690 srgb} bind def
/col23 {0.820 0.000 0.820 srgb} bind def
/col24 {0.500 0.190 0.000 srgb} bind def
/col25 {0.630 0.250 0.000 srgb} bind def
/col26 {0.750 0.380 0.000 srgb} bind def
/col27 {1.000 0.500 0.500 srgb} bind def
/col28 {1.000 0.630 0.630 srgb} bind def
/col29 {1.000 0.750 0.750 srgb} bind def
/col30 {1.000 0.880 0.880 srgb} bind def
/col31 {1.000 0.840 0.000 srgb} bind def

end
save
-3.0 459.0 translate
1 -1 scale

/cp {closepath} bind def
/ef {eofill} bind def
/gr {grestore} bind def
/gs {gsave} bind def
/sa {save} bind def
/rs {restore} bind def
/l {lineto} bind def
/m {moveto} bind def
/rm {rmoveto} bind def
/n {newpath} bind def
/s {stroke} bind def
/sh {show} bind def
/slc {setlinecap} bind def
/slj {setlinejoin} bind def
/slw {setlinewidth} bind def
/srgb {setrgbcolor} bind def
/rot {rotate} bind def
/sc {scale} bind def
/sd {setdash} bind def
/ff {findfont} bind def
/sf {setfont} bind def
/scf {scalefont} bind def
/sw {stringwidth} bind def
/tr {translate} bind def
/tnt {dup dup currentrgbcolor
  4 -2 roll dup 1 exch sub 3 -1 roll mul add
  4 -2 roll dup 1 exch sub 3 -1 roll mul add
  4 -2 roll dup 1 exch sub 3 -1 roll mul add srgb}
  bind def
/shd {dup dup currentrgbcolor 4 -2 roll mul 4 -2 roll mul
  4 -2 roll mul srgb} bind def
 /DrawEllipse {
/endangle exch def
/startangle exch def
/yrad exch def
/xrad exch def
/y exch def
/x exch def
/savematrix mtrx currentmatrix def
x y tr xrad yrad sc 0 0 1 startangle endangle arc
closepath
savematrix setmatrix
} def

/$F2psBegin {$F2psDict begin /$F2psEnteredState save def} def
/$F2psEnd {$F2psEnteredState restore end} def

$F2psBegin
10 setmiterlimit
n 0 792 m 0 0 l 612 0 l 612 792 l cp clip
 0.06000 0.06000 sc
30.000 slw
n 6980 1150 255 255 0 360 DrawEllipse gs col-1 s gr

n 8383 1320 341 341 0 360 DrawEllipse gs col-1 s gr

n 7093 2540 216 216 0 360 DrawEllipse gs col-1 s gr

n 8610 2440 323 323 0 360 DrawEllipse gs col-1 s gr

7.500 slw
n 6931 1399 40 40 0 360 DrawEllipse gs col7 0.00 shd ef gr gs col-1 s gr

n 6998 2347 40 40 0 360 DrawEllipse gs col7 0.00 shd ef gr gs col-1 s gr

n 8559 1612 40 40 0 360 DrawEllipse gs col7 0.00 shd ef gr gs col-1 s gr

n 8658 2122 40 40 0 360 DrawEllipse gs col7 0.00 shd ef gr gs col-1 s gr

30.000 slw
gs n 8114 87 m
8501.6 159.7 9253.9 462.4 9532 796 curveto
9842.9 1168.1 10216.1 1947.1 10028 2497 curveto
9823.9 3095.7 8980.2 3514.7 8469 3632 curveto
7837.6 3775.6 6676.6 3702.4 6129 3206 curveto
5653.6 2775.4 5404.5 1828.1 5562 1221 curveto
5668.8 808.2 6223.0 383.6 6554 228 curveto
6924.1 55.2 7731.2 14.4 8114 87 curveto
 cp  gs col-1 s gr
 gr

1 slc
15.000 slw
gs n 6930 1399 m
6874.5 1712.9 6859.5 1850.7 6870 1950 curveto
6877.9 2024.9 6910.7 2124.2 7001 2347 curveto
 gs col-1 s gr
 gr

gs n 8560 1612 m
8659.5 1744.9 8696.0 1807.9 8706 1864 curveto
8715.0 1914.5 8703.7 1979.0 8661 2122 curveto
 gs col-1 s gr
 gr

30.000 slw
n 8383 5219 341 341 0 360 DrawEllipse gs col-1 s gr

n 7093 6440 216 216 0 360 DrawEllipse gs col-1 s gr

n 8610 6340 323 323 0 360 DrawEllipse gs col-1 s gr

n 6980 5049 255 255 0 360 DrawEllipse gs col-1 s gr

7.500 slw
n 6846 4840 40 40 0 360 DrawEllipse gs col7 0.00 shd ef gr gs col-1 s gr

n 6512 4150 40 40 0 360 DrawEllipse gs col7 0.00 shd ef gr gs col-1 s gr

n 5537 5790 40 40 0 360 DrawEllipse gs col7 0.00 shd ef gr gs col-1 s gr

n 7299 6367 40 40 0 360 DrawEllipse gs col7 0.00 shd ef gr gs col-1 s gr

n 8139 5463 40 40 0 360 DrawEllipse gs col7 0.00 shd ef gr gs col-1 s gr

n 8326 6176 40 40 0 360 DrawEllipse gs col7 0.00 shd ef gr gs col-1 s gr

30.000 slw
gs n 8114 3986 m
8501.4 4058.8 9253.7 4361.5 9532 4695 curveto
9842.7 5067.4 10215.9 5846.7 10028 6397 curveto
9823.6 6995.5 8980.0 7414.7 8468 7532 curveto
7837.4 7675.6 6676.4 7602.4 6129 7106 curveto
5653.4 6675.1 5404.3 5727.6 5562 5120 curveto
5668.6 4707.4 6222.8 4282.7 6554 4127 curveto
6923.9 3954.2 7731.0 3913.4 8114 3986 curveto
 cp  gs col-1 s gr
 gr

15.000 slw
gs n 8327 6180 m
7895.1 5865.2 7693.1 5747.6 7519 5709 curveto
7271.5 5655.2 6755.8 5801.1 6527 5805 curveto
6394.9 5806.4 6084.6 5767.2 5953 5766 curveto
5881.5 5766.1 5777.5 5771.9 5537 5790 curveto
 gs col-1 s gr
 gr

gs n 8132 5459 m
8008.2 5580.2 7956.2 5634.0 7924 5674 curveto
7829.5 5793.5 7685.1 6119.0 7566 6229 curveto
7525.2 6266.2 7458.7 6301.7 7300 6371 curveto
 gs col-1 s gr
 gr

gs n 6838 4836 m
6740.8 4652.9 6699.5 4573.4 6673 4518 curveto
6643.5 4456.8 6603.3 4365.7 6512 4154 curveto
 gs col-1 s gr
 gr

30.000 slw
n 1578 5050 255 255 0 360 DrawEllipse gs col-1 s gr

n 2981 5220 341 341 0 360 DrawEllipse gs col-1 s gr

n 1691 6439 216 216 0 360 DrawEllipse gs col-1 s gr

n 3208 6340 323 323 0 360 DrawEllipse gs col-1 s gr

7.500 slw
n 2973 5559 40 40 0 360 DrawEllipse gs col7 0.00 shd ef gr gs col-1 s gr

n 3045 6062 40 40 0 360 DrawEllipse gs col7 0.00 shd ef gr gs col-1 s gr

n 2902 6459 40 40 0 360 DrawEllipse gs col7 0.00 shd ef gr gs col-1 s gr

n 1890 6523 40 40 0 360 DrawEllipse gs col7 0.00 shd ef gr gs col-1 s gr

n 3510 6457 40 40 0 360 DrawEllipse gs col7 0.00 shd ef gr gs col-1 s gr

n 4257 6914 40 40 0 360 DrawEllipse gs col7 0.00 shd ef gr gs col-1 s gr

n 1612 4005 40 40 0 360 DrawEllipse gs col7 0.00 shd ef gr gs col-1 s gr

n 1432 4841 40 40 0 360 DrawEllipse gs col7 0.00 shd ef gr gs col-1 s gr

30.000 slw
gs n 2712 3986 m
3099.4 4059.2 3851.7 4361.8 4130 4695 curveto
4440.7 5067.6 4813.9 5846.6 4626 6397 curveto
4421.6 6995.1 3578.0 7414.2 3066 7531 curveto
2435.4 7675.1 1274.4 7601.8 727 7106 curveto
251.4 6674.8 2.3 5727.6 160 5120 curveto
266.6 4707.6 820.8 4283.1 1152 4128 curveto
1521.8 3954.7 2329.0 3913.8 2712 3986 curveto
 cp  gs col-1 s gr
 gr

15.000 slw
gs n 3505 6455 m
3790.5 6557.7 3911.8 6608.9 3990 6660 curveto
4045.2 6696.0 4112.7 6759.8 4260 6915 curveto
 gs col-1 s gr
 gr

gs n 2972 5563 m
2973.0 5737.6 2976.2 5813.1 2985 5865 curveto
2990.9 5900.1 3005.7 5949.1 3044 6061 curveto
 gs col-1 s gr
 gr

gs n 1889 6526 m
2199.1 6624.3 2338.1 6654.1 2445 6645 curveto
2536.2 6637.2 2648.7 6592.2 2895 6465 curveto
 gs col-1 s gr
 gr

gs n 1429 4838 m
1311.7 4643.5 1293.7 4549.7 1357 4463 curveto
1407.1 4394.4 1526.7 4494.4 1582 4436 curveto
1659.5 4354.0 1666.3 4248.3 1609 4013 curveto
 gs col-1 s gr
 gr

30.000 slw
n 2280 1860 511 511 0 360 DrawEllipse gs col-1 s gr

7.500 slw
n 1270 185 40 40 0 360 DrawEllipse gs col7 0.00 shd ef gr gs col-1 s gr

n 1770 1810 40 40 0 360 DrawEllipse gs col7 0.00 shd ef gr gs col-1 s gr

15.000 slw
gs n 1265 175 m
1492.2 561.5 1567.2 744.0 1565 905 curveto
1562.3 1102.9 1173.7 1308.5 1270 1570 curveto
1325.1 1719.5 1450.1 1779.5 1770 1810 curveto
 gs col-1 s gr
 gr

30.000 slw
gs n 2712 87 m
3099.6 160.2 3851.9 462.7 4130 796 curveto
4440.9 1168.3 4814.1 1947.1 4626 2498 curveto
4421.9 3095.3 3578.2 3514.2 3067 3631 curveto
2435.6 3775.1 1274.6 3701.8 727 3206 curveto
251.6 2775.1 2.5 1828.1 160 1221 curveto
266.8 808.4 821.0 383.9 1152 229 curveto
1522.1 55.7 2329.2 14.8 2712 87 curveto
 cp  gs col-1 s gr
 gr

$F2psEnd
rs
EndOfTheIncludedPostscriptMagicCookie

\closepsdump

\psdump{evenodd.eps}%!PS-Adobe-2.0 EPSF-2.0
/$F2psDict 200 dict def
$F2psDict begin
$F2psDict /mtrx matrix put
/col-1 {0 setgray} bind def
/col0 {0.000 0.000 0.000 srgb} bind def
/col1 {0.000 0.000 1.000 srgb} bind def
/col2 {0.000 1.000 0.000 srgb} bind def
/col3 {0.000 1.000 1.000 srgb} bind def
/col4 {1.000 0.000 0.000 srgb} bind def
/col5 {1.000 0.000 1.000 srgb} bind def
/col6 {1.000 1.000 0.000 srgb} bind def
/col7 {1.000 1.000 1.000 srgb} bind def
/col8 {0.000 0.000 0.560 srgb} bind def
/col9 {0.000 0.000 0.690 srgb} bind def
/col10 {0.000 0.000 0.820 srgb} bind def
/col11 {0.530 0.810 1.000 srgb} bind def
/col12 {0.000 0.560 0.000 srgb} bind def
/col13 {0.000 0.690 0.000 srgb} bind def
/col14 {0.000 0.820 0.000 srgb} bind def
/col15 {0.000 0.560 0.560 srgb} bind def
/col16 {0.000 0.690 0.690 srgb} bind def
/col17 {0.000 0.820 0.820 srgb} bind def
/col18 {0.560 0.000 0.000 srgb} bind def
/col19 {0.690 0.000 0.000 srgb} bind def
/col20 {0.820 0.000 0.000 srgb} bind def
/col21 {0.560 0.000 0.560 srgb} bind def
/col22 {0.690 0.000 0.690 srgb} bind def
/col23 {0.820 0.000 0.820 srgb} bind def
/col24 {0.500 0.190 0.000 srgb} bind def
/col25 {0.630 0.250 0.000 srgb} bind def
/col26 {0.750 0.380 0.000 srgb} bind def
/col27 {1.000 0.500 0.500 srgb} bind def
/col28 {1.000 0.630 0.630 srgb} bind def
/col29 {1.000 0.750 0.750 srgb} bind def
/col30 {1.000 0.880 0.880 srgb} bind def
/col31 {1.000 0.840 0.000 srgb} bind def

end
save
-33.0 213.0 translate
1 -1 scale

/cp {closepath} bind def
/ef {eofill} bind def
/gr {grestore} bind def
/gs {gsave} bind def
/sa {save} bind def
/rs {restore} bind def
/l {lineto} bind def
/m {moveto} bind def
/rm {rmoveto} bind def
/n {newpath} bind def
/s {stroke} bind def
/sh {show} bind def
/slc {setlinecap} bind def
/slj {setlinejoin} bind def
/slw {setlinewidth} bind def
/srgb {setrgbcolor} bind def
/rot {rotate} bind def
/sc {scale} bind def
/sd {setdash} bind def
/ff {findfont} bind def
/sf {setfont} bind def
/scf {scalefont} bind def
/sw {stringwidth} bind def
/tr {translate} bind def
/tnt {dup dup currentrgbcolor
  4 -2 roll dup 1 exch sub 3 -1 roll mul add
  4 -2 roll dup 1 exch sub 3 -1 roll mul add
  4 -2 roll dup 1 exch sub 3 -1 roll mul add srgb}
  bind def
/shd {dup dup currentrgbcolor 4 -2 roll mul 4 -2 roll mul
  4 -2 roll mul srgb} bind def
 /DrawEllipse {
/endangle exch def
/startangle exch def
/yrad exch def
/xrad exch def
/y exch def
/x exch def
/savematrix mtrx currentmatrix def
x y tr xrad yrad sc 0 0 1 startangle endangle arc
closepath
savematrix setmatrix
} def

/$F2psBegin {$F2psDict begin /$F2psEnteredState save def} def
/$F2psEnd {$F2psEnteredState restore end} def

$F2psBegin
10 setmiterlimit
n 0 792 m 0 0 l 612 0 l 612 792 l cp clip
 0.06000 0.06000 sc
30.000 slw
n 6843 1741 201 201 0 360 DrawEllipse gs col-1 s gr

n 7946 1875 268 268 0 360 DrawEllipse gs col-1 s gr

n 6932 2834 170 170 0 360 DrawEllipse gs col-1 s gr

7.500 slw
n 7906 1615 34 34 0 360 DrawEllipse gs col7 0.00 shd ef gr gs col-1 s gr

n 7453 1039 34 34 0 360 DrawEllipse gs col7 0.00 shd ef gr gs col-1 s gr

n 8510 2958 34 34 0 360 DrawEllipse gs col7 0.00 shd ef gr gs col-1 s gr

30.000 slw
n 2855 1626 268 268 0 360 DrawEllipse gs col-1 s gr

n 1841 2584 170 170 0 360 DrawEllipse gs col-1 s gr

n 3033 2506 254 254 0 360 DrawEllipse gs col-1 s gr

7.500 slw
n 3140 822 34 34 0 360 DrawEllipse gs col7 0.00 shd ef gr gs col-1 s gr

n 2818 1366 34 34 0 360 DrawEllipse gs col7 0.00 shd ef gr gs col-1 s gr

30.000 slw
n 1751 1492 201 201 0 360 DrawEllipse gs col-1 s gr

7.500 slw
n 1865 661 34 34 0 360 DrawEllipse gs col7 0.00 shd ef gr gs col-1 s gr

n 1756 2435 34 34 0 360 DrawEllipse gs col7 0.00 shd ef gr gs col-1 s gr

n 2788 2427 34 34 0 360 DrawEllipse gs col7 0.00 shd ef gr gs col-1 s gr

n 1946 1470 34 34 0 360 DrawEllipse gs col7 0.00 shd ef gr gs col-1 s gr

30.000 slw
n 7350 2025 1489 1489 0 360 DrawEllipse gs col-1 s gr

7.500 slw
n 6845 2683 34 34 0 360 DrawEllipse gs col7 0.00 shd ef gr gs col-1 s gr

n 7035 1718 34 34 0 360 DrawEllipse gs col7 0.00 shd ef gr gs col-1 s gr

1 slj
1 slc
gs  clippath
5250 2033 m 5370 2063 l 5250 2093 l 5412 2093 l 5412 2033 l  cp clip
n 4505 2063 m 5397 2063 l gs col7 0.00 shd ef gr gs col-1 s gr gr

0 slc
0 slj
n 5250 2033 m 5370 2063 l 5250 2093 l 5250 2063 l 5250 2033 l  cp gs 0.00 setgray ef gr  col-1 s
30.000 slw
gs n 2643 656 m
2948.0 713.2 3539.5 951.1 3758 1213 curveto
4002.6 1506.0 4296.0 2118.4 4148 2551 curveto
3987.6 3021.2 3324.3 3350.7 2922 3443 curveto
2425.9 3555.8 1513.0 3498.2 1082 3108 curveto
708.6 2769.4 512.8 2024.8 636 1548 curveto
720.6 1223.0 1156.4 889.2 1417 767 curveto
1707.6 631.1 2342.2 598.9 2643 656 curveto
 cp  gs col-1 s gr
 gr

1 slc
15.000 slw
gs n 7907 1613 m
7899.3 1528.9 7892.4 1492.7 7879 1468 curveto
7837.8 1389.2 7686.1 1275.0 7630 1218 curveto
7598.9 1186.3 7554.4 1140.2 7451 1033 curveto
 gs col-1 s gr
 gr

gs n 2816 1363 m
2785.6 1285.3 2778.7 1249.1 2788 1218 curveto
2824.0 1102.0 3023.4 999.1 3083 923 curveto
3095.9 907.2 3111.3 882.1 3145 823 curveto
 gs col-1 s gr
 gr

7.500 slw
gs  clippath
7226 900 m 7333 962 l 7209 958 l 7365 1003 l 7382 945 l  cp clip
 [66.7] 0 sd
n 3925 1260 m
4616.5 952.8 4923.3 835.7 5152 792 curveto
5477.5 729.0 6231.9 732.1 6556 769 curveto
6699.7 785.9 6900.4 836.1 7359 970 curveto
 gs col-1 s gr
 gr
 [] 0 sd
0 slc
n 7226 900 m 7333 962 l 7209 958 l 7218 929 l 7226 900 l  cp gs 0.00 setgray ef gr  col-1 s
1 slc
gs  clippath
5954 3040 m 6066 2987 l 5992 3086 l 6118 2984 l 6080 2937 l  cp clip
 [66.7] 0 sd
n 3234 2765 m
3744.4 3185.7 3984.2 3347.2 4193 3411 curveto
4490.9 3502.0 5132.2 3501.7 5430 3411 curveto
5572.9 3367.5 5737.1 3257.2 6087 2970 curveto
 gs col-1 s gr
 gr
 [] 0 sd
0 slc
n 5954 3040 m 6066 2987 l 5992 3086 l 5973 3063 l 5954 3040 l  cp gs 0.00 setgray ef gr  col-1 s
1 slc
15.000 slw
gs n 1950 1470 m
1981.6 1409.8 1992.8 1382.2 1995 1360 curveto
2009.1 1216.3 1907.8 924.2 1885 795 curveto
1880.8 771.3 1875.3 736.6 1863 656 curveto
 gs col-1 s gr
 gr

gs n 1945 1475 m
1987.7 1525.5 2002.7 1550.5 2005 1575 curveto
2027.0 1808.1 1705.4 2128.5 1718 2355 curveto
1719.0 2372.8 1728.7 2392.3 1757 2433 curveto
 gs col-1 s gr
 gr

gs n 1945 1470 m
2011.7 1457.4 2041.7 1457.4 2065 1470 curveto
2351.0 1624.5 2436.8 2223.2 2671 2395 curveto
2690.3 2409.2 2719.6 2417.4 2788 2428 curveto
 gs col-1 s gr
 gr

gs n 7039 1718 m
7064.2 1654.0 7075.5 1626.5 7084 1608 curveto
7137.6 1491.8 7248.6 1212.3 7341 1110 curveto
7358.0 1091.2 7386.7 1072.5 7456 1035 curveto
 gs col-1 s gr
 gr

gs n 7034 1723 m
7076.7 1773.5 7091.7 1798.5 7094 1823 curveto
7116.0 2056.1 6794.4 2376.5 6807 2603 curveto
6808.0 2620.8 6817.7 2640.3 6846 2681 curveto
 gs col-1 s gr
 gr

gs n 7034 1718 m
7100.7 1705.1 7130.7 1705.1 7154 1718 curveto
7424.5 1867.1 7508.4 2429.7 7711 2600 curveto
7853.9 2720.1 8240.7 2802.8 8391 2880 curveto
8412.7 2891.1 8442.7 2909.9 8511 2955 curveto
 gs col-1 s gr
 gr

$F2psEnd
rs
EndOfTheIncludedPostscriptMagicCookie

\closepsdump

\psdump{furtherx.eps}%!PS-Adobe-2.0 EPSF-2.0
/$F2psDict 200 dict def
$F2psDict begin
$F2psDict /mtrx matrix put
/col-1 {0 setgray} bind def
/col0 {0.000 0.000 0.000 srgb} bind def
/col1 {0.000 0.000 1.000 srgb} bind def
/col2 {0.000 1.000 0.000 srgb} bind def
/col3 {0.000 1.000 1.000 srgb} bind def
/col4 {1.000 0.000 0.000 srgb} bind def
/col5 {1.000 0.000 1.000 srgb} bind def
/col6 {1.000 1.000 0.000 srgb} bind def
/col7 {1.000 1.000 1.000 srgb} bind def
/col8 {0.000 0.000 0.560 srgb} bind def
/col9 {0.000 0.000 0.690 srgb} bind def
/col10 {0.000 0.000 0.820 srgb} bind def
/col11 {0.530 0.810 1.000 srgb} bind def
/col12 {0.000 0.560 0.000 srgb} bind def
/col13 {0.000 0.690 0.000 srgb} bind def
/col14 {0.000 0.820 0.000 srgb} bind def
/col15 {0.000 0.560 0.560 srgb} bind def
/col16 {0.000 0.690 0.690 srgb} bind def
/col17 {0.000 0.820 0.820 srgb} bind def
/col18 {0.560 0.000 0.000 srgb} bind def
/col19 {0.690 0.000 0.000 srgb} bind def
/col20 {0.820 0.000 0.000 srgb} bind def
/col21 {0.560 0.000 0.560 srgb} bind def
/col22 {0.690 0.000 0.690 srgb} bind def
/col23 {0.820 0.000 0.820 srgb} bind def
/col24 {0.500 0.190 0.000 srgb} bind def
/col25 {0.630 0.250 0.000 srgb} bind def
/col26 {0.750 0.380 0.000 srgb} bind def
/col27 {1.000 0.500 0.500 srgb} bind def
/col28 {1.000 0.630 0.630 srgb} bind def
/col29 {1.000 0.750 0.750 srgb} bind def
/col30 {1.000 0.880 0.880 srgb} bind def
/col31 {1.000 0.840 0.000 srgb} bind def

end
save
-34.0 273.0 translate
1 -1 scale

/cp {closepath} bind def
/ef {eofill} bind def
/gr {grestore} bind def
/gs {gsave} bind def
/sa {save} bind def
/rs {restore} bind def
/l {lineto} bind def
/m {moveto} bind def
/rm {rmoveto} bind def
/n {newpath} bind def
/s {stroke} bind def
/sh {show} bind def
/slc {setlinecap} bind def
/slj {setlinejoin} bind def
/slw {setlinewidth} bind def
/srgb {setrgbcolor} bind def
/rot {rotate} bind def
/sc {scale} bind def
/sd {setdash} bind def
/ff {findfont} bind def
/sf {setfont} bind def
/scf {scalefont} bind def
/sw {stringwidth} bind def
/tr {translate} bind def
/tnt {dup dup currentrgbcolor
  4 -2 roll dup 1 exch sub 3 -1 roll mul add
  4 -2 roll dup 1 exch sub 3 -1 roll mul add
  4 -2 roll dup 1 exch sub 3 -1 roll mul add srgb}
  bind def
/shd {dup dup currentrgbcolor 4 -2 roll mul 4 -2 roll mul
  4 -2 roll mul srgb} bind def
 /DrawEllipse {
/endangle exch def
/startangle exch def
/yrad exch def
/xrad exch def
/y exch def
/x exch def
/savematrix mtrx currentmatrix def
x y tr xrad yrad sc 0 0 1 startangle endangle arc
closepath
savematrix setmatrix
} def

/$F2psBegin {$F2psDict begin /$F2psEnteredState save def} def
/$F2psEnd {$F2psEnteredState restore end} def

$F2psBegin
10 setmiterlimit
n 0 792 m 0 0 l 612 0 l 612 792 l cp clip
 0.06000 0.06000 sc
30.000 slw
n 2295 3270 229 229 0 360 DrawEllipse gs col-1 s gr

n 3900 3165 342 342 0 360 DrawEllipse gs col-1 s gr

n 3000 1800 361 361 0 360 DrawEllipse gs col-1 s gr

7.500 slw
n 3075 2152 54 54 0 360 DrawEllipse gs col7 0.00 shd ef gr gs col-1 s gr

n 2595 4477 54 54 0 360 DrawEllipse gs col7 0.00 shd ef gr gs col-1 s gr

n 3270 2043 54 54 0 360 DrawEllipse gs col7 0.00 shd ef gr gs col-1 s gr

n 3810 4413 54 54 0 360 DrawEllipse gs col7 0.00 shd ef gr gs col-1 s gr

1 slj
1 slc
15.000 slw
n 2922 1357 m 3000 1433 l 2925 1508 l gs col-1 s gr
n 2892 867 m 2970 943 l 2895 1018 l gs col-1 s gr
n 2459 2237 m 2410 2333 l 2316 2285 l gs col-1 s gr
n 4976 3693 m 5083 3703 l 5074 3807 l gs col-1 s gr
n 2911 3178 m 3007 3129 l 3055 3223 l gs col-1 s gr
n 3041 713 m 2963 637 l 3038 562 l gs col-1 s gr
n 1015 3499 m 980 3601 l 880 3567 l gs col-1 s gr
n 4824 2924 m 4735 2987 l 4674 2900 l gs col-1 s gr
30.000 slw
gs n 3375 675 m
3785.0 752.4 4580.8 1072.5 4875 1425 curveto
5203.9 1819.1 5598.7 2643.2 5400 3225 curveto
5183.8 3858.1 4291.3 4301.4 3750 4425 curveto
3082.5 4577.4 1854.3 4499.9 1275 3975 curveto
772.1 3519.3 508.6 2517.3 675 1875 curveto
788.1 1438.3 1374.5 989.2 1725 825 curveto
2116.1 641.8 2969.9 598.6 3375 675 curveto
 cp  gs col-1 s gr
 gr

15.000 slw
gs n 3075 2145 m
3126.2 2392.3 3141.2 2501.0 3135 2580 curveto
3101.4 3010.0 2705.4 3895.2 2617 4297 curveto
2610.7 4325.9 2604.7 4368.4 2593 4467 curveto
 gs col-1 s gr
 gr

0 slc
gs n 3265 2038 m
3322.8 2115.5 3357.8 2138.0 3405 2128 curveto
3631.6 2080.2 3798.6 1823.1 3735 1628 curveto
3621.1 1278.7 3107.7 1084.4 2730 1253 curveto
2378.2 1410.0 2351.3 1907.2 2365 2183 curveto
2379.4 2472.7 2867.9 2793.8 2843 3117 curveto
2826.3 3333.8 2682.7 3631.8 2430 3687 curveto
2090.1 3761.3 1759.7 3434.4 1605 3222 curveto
1440.0 2995.4 1361.9 2494.0 1403 2232 curveto
1439.1 2002.0 1681.0 1614.3 1838 1459 curveto
2022.7 1276.3 2471.1 985.5 2753 957 curveto
3211.1 910.7 4026.1 1079.1 4373 1459 curveto
4711.8 1830.0 4826.1 2626.0 4718 3079 curveto
4620.5 3487.8 3888.7 3849.4 3780 4243 curveto
3771.0 4275.7 3776.7 4316.7 3803 4407 curveto
 gs col-1 s gr
 gr

$F2psEnd
rs
EndOfTheIncludedPostscriptMagicCookie

\closepsdump

\psdump{kevelines.eps}%!PS-Adobe-2.0 EPSF-2.0
/$F2psDict 200 dict def
$F2psDict begin
$F2psDict /mtrx matrix put
/col-1 {0 setgray} bind def
/col0 {0.000 0.000 0.000 srgb} bind def
/col1 {0.000 0.000 1.000 srgb} bind def
/col2 {0.000 1.000 0.000 srgb} bind def
/col3 {0.000 1.000 1.000 srgb} bind def
/col4 {1.000 0.000 0.000 srgb} bind def
/col5 {1.000 0.000 1.000 srgb} bind def
/col6 {1.000 1.000 0.000 srgb} bind def
/col7 {1.000 1.000 1.000 srgb} bind def
/col8 {0.000 0.000 0.560 srgb} bind def
/col9 {0.000 0.000 0.690 srgb} bind def
/col10 {0.000 0.000 0.820 srgb} bind def
/col11 {0.530 0.810 1.000 srgb} bind def
/col12 {0.000 0.560 0.000 srgb} bind def
/col13 {0.000 0.690 0.000 srgb} bind def
/col14 {0.000 0.820 0.000 srgb} bind def
/col15 {0.000 0.560 0.560 srgb} bind def
/col16 {0.000 0.690 0.690 srgb} bind def
/col17 {0.000 0.820 0.820 srgb} bind def
/col18 {0.560 0.000 0.000 srgb} bind def
/col19 {0.690 0.000 0.000 srgb} bind def
/col20 {0.820 0.000 0.000 srgb} bind def
/col21 {0.560 0.000 0.560 srgb} bind def
/col22 {0.690 0.000 0.690 srgb} bind def
/col23 {0.820 0.000 0.820 srgb} bind def
/col24 {0.500 0.190 0.000 srgb} bind def
/col25 {0.630 0.250 0.000 srgb} bind def
/col26 {0.750 0.380 0.000 srgb} bind def
/col27 {1.000 0.500 0.500 srgb} bind def
/col28 {1.000 0.630 0.630 srgb} bind def
/col29 {1.000 0.750 0.750 srgb} bind def
/col30 {1.000 0.880 0.880 srgb} bind def
/col31 {1.000 0.840 0.000 srgb} bind def

end
save
-2.0 227.0 translate
1 -1 scale

/cp {closepath} bind def
/ef {eofill} bind def
/gr {grestore} bind def
/gs {gsave} bind def
/sa {save} bind def
/rs {restore} bind def
/l {lineto} bind def
/m {moveto} bind def
/rm {rmoveto} bind def
/n {newpath} bind def
/s {stroke} bind def
/sh {show} bind def
/slc {setlinecap} bind def
/slj {setlinejoin} bind def
/slw {setlinewidth} bind def
/srgb {setrgbcolor} bind def
/rot {rotate} bind def
/sc {scale} bind def
/sd {setdash} bind def
/ff {findfont} bind def
/sf {setfont} bind def
/scf {scalefont} bind def
/sw {stringwidth} bind def
/tr {translate} bind def
/tnt {dup dup currentrgbcolor
  4 -2 roll dup 1 exch sub 3 -1 roll mul add
  4 -2 roll dup 1 exch sub 3 -1 roll mul add
  4 -2 roll dup 1 exch sub 3 -1 roll mul add srgb}
  bind def
/shd {dup dup currentrgbcolor 4 -2 roll mul 4 -2 roll mul
  4 -2 roll mul srgb} bind def
 /DrawEllipse {
/endangle exch def
/startangle exch def
/yrad exch def
/xrad exch def
/y exch def
/x exch def
/savematrix mtrx currentmatrix def
x y tr xrad yrad sc 0 0 1 startangle endangle arc
closepath
savematrix setmatrix
} def

/$F2psBegin {$F2psDict begin /$F2psEnteredState save def} def
/$F2psEnd {$F2psEnteredState restore end} def

$F2psBegin
10 setmiterlimit
n 0 792 m 0 0 l 612 0 l 612 792 l cp clip
 0.06000 0.06000 sc
30.000 slw
n 1568 1179 255 255 0 360 DrawEllipse gs col-1 s gr

n 2971 1349 341 341 0 360 DrawEllipse gs col-1 s gr

n 1681 2569 216 216 0 360 DrawEllipse gs col-1 s gr

n 3198 2469 323 323 0 360 DrawEllipse gs col-1 s gr

7.500 slw
n 2870 2440 40 40 0 360 DrawEllipse gs col7 0.00 shd ef gr gs col-1 s gr

n 1600 1440 40 40 0 360 DrawEllipse gs col7 0.00 shd ef gr gs col-1 s gr

n 1665 2355 40 40 0 360 DrawEllipse gs col7 0.00 shd ef gr gs col-1 s gr

n 3160 250 40 40 0 360 DrawEllipse gs col7 0.00 shd ef gr gs col-1 s gr

n 3005 1010 40 40 0 360 DrawEllipse gs col7 0.00 shd ef gr gs col-1 s gr

n 535 3030 40 40 0 360 DrawEllipse gs col7 0.00 shd ef gr gs col-1 s gr

30.000 slw
gs n 2702 116 m
3089.6 188.7 3841.9 491.4 4120 825 curveto
4430.9 1197.1 4804.1 1976.1 4616 2526 curveto
4411.9 3124.7 3568.2 3543.7 3057 3661 curveto
2425.6 3804.6 1264.6 3731.4 717 3235 curveto
241.6 2804.4 -7.5 1857.1 150 1250 curveto
256.8 837.2 811.0 412.6 1142 257 curveto
1512.1 84.2 2319.2 43.4 2702 116 curveto
 cp  gs col-1 s gr
 gr

1 slc
15.000 slw
gs n 1610 1445 m
1674.8 1474.1 1702.3 1487.8 1720 1500 curveto
1976.3 1676.1 2417.5 2245.7 2705 2400 curveto
2731.7 2414.3 2771.7 2424.3 2865 2440 curveto
 gs col-1 s gr
 gr

gs n 3150 255 m
3121.7 305.7 3110.5 328.2 3105 345 curveto
3063.8 471.2 3024.0 775.8 3010 905 curveto
3008.1 922.9 3006.8 949.2 3005 1010 curveto
 gs col-1 s gr
 gr

gs n 1590 1440 m
1527.7 1487.0 1501.4 1508.2 1485 1525 curveto
1350.0 1663.1 1072.5 2017.3 975 2185 curveto
882.5 2344.0 777.0 2758.8 670 2915 curveto
651.7 2941.7 619.2 2973.0 540 3040 curveto
 gs col-1 s gr
 gr

gs n 1600 1440 m
1598.8 1518.1 1598.8 1551.9 1600 1575 curveto
1608.0 1729.9 1658.2 2095.7 1670 2250 curveto
1671.3 2267.1 1672.6 2292.1 1675 2350 curveto
 gs col-1 s gr
 gr

$F2psEnd
rs
EndOfTheIncludedPostscriptMagicCookie

\closepsdump

\newcommand\ignore[1]{}
\usepackage[dvips]{epsfig}
\usepackage{amstex,ifthen,amssymb,alltt,calc,epic,theorem}
\theoremheaderfont\scshape
\DeclareFontFamily{OT1}{rsfs}{}
\DeclareFontShape{OT1}{rsfs}{m}{n}{ <-7> rsfs5 <7-10> rsfs7 <10-> rsfs10}{}
\DeclareMathAlphabet\mathcurl{OT1}{rsfs}{m}{n}

\newcommand\afigwidth0%
\newcounter{afigheight}%
\newcommand\bfigwidth0%
\newcounter{bfigheight}%
\newcommand\aboundingbox[4]{%
  \setcounter{afigheight}{\afigwidth*(#4-#2)/(#3-#1)}}
\newcommand\bboundingbox[4]{%
  \setcounter{bfigheight}{\afigwidth*(#4-#2)/(#3-#1)}}
\newlength\aeqspace%
\newlength\beqspace%
\newlength\border%
\newlength\tempeqspace%
\newcounter{tempactr}%
\newcounter{tempbctr}%
\newcounter{tempcctr}%
2%
\newcommand{\fig}[6]{%
  \setcounter{tempactr}{#6}%
  \setcounter{tempbctr}{200+\border/\unitlength}%
  \setcounter{tempcctr}{#6+\border/\unitlength}%
  \setlength\tempeqspace{#5}%
  \vskip\tempeqspace%
  \vskip\border%
  \centering%
  \begin{picture}(200,\value{tempactr})%
    \put(0,0){\epsfig{figure=#1,%
                      width=200\unitlength}}%
    #3%
  \end{picture}%
  \vskip\border%
  \vskip\tempeqspace%
  \vskip\abovecaptionskip%
  \refstepcounter{figure}#2%
  \makebox[0pt][c]{\parbox[t]{1.5\textwidth}{\centering%
    \ifthenelse{\equal{#4}{}}{Figure \thefigure}{Figure \thefigure: #4}}}%
  \vskip\belowcaptionskip}%

\newcommand{\singdiag}[2]{{%
  \setlength\border{0pt}%
  \renewcommand\baselinestretch1\normalsize%
  #1%
  \begin{figure}[t]%
    \par%
    \noindent%
    \hfill%
    \setlength\unitlength{\textwidth/100*\afigwidth/200}%
    \begin{minipage}[t]{\textwidth/100*\afigwidth+2\border}%
      #2{0pt}{200*\value{afigheight}/\afigwidth}%
    \end{minipage}%
    \hfill%
    \null%
  \end{figure}}}%

\newcommand{\doubdiag}[3]{{\setlength\border{0pt}%{\textwidth/100*2}%
  \renewcommand\baselinestretch1\normalsize%
  \setlength\aeqspace{0pt}%
  \setlength\beqspace{0pt}%
  #1%
  \ifthenelse{\value{bfigheight} > \value{afigheight}}{%
    \setlength\aeqspace{\textwidth/200*(\value{bfigheight}-%
      \value{afigheight})}}{%
    \setlength\beqspace{\textwidth/200*(\value{afigheight}-%
      \value{bfigheight})}}%
  \begin{figure}[t]%
    \par%
    \noindent%
    \hfill%
    \setlength\unitlength{\textwidth/100*\afigwidth/200}%
    \begin{minipage}[t]{\textwidth/100*\afigwidth+2\border}%
      #2\aeqspace{200*\value{afigheight}/\afigwidth}%
    \end{minipage}%
    \hfill%
    \setlength\unitlength{\textwidth/100*\bfigwidth/200}%
    \begin{minipage}[t]{\textwidth/100*\bfigwidth+2\border}%
      #3\beqspace{200*\value{bfigheight}/\bfigwidth}%
    \end{minipage}%
    \hfill%
    \null%
  \end{figure}}}%

\newcommand\x{\mathbf x}
\renewcommand\r{\mathbb R}
\renewcommand\c{\mathbb C}
\newcommand\z{\mathbb Z}

\renewcommand\d{\mathrm d}
\newcommand\symm{G}
\newcommand\grad\nabla

\newcommand\qed{\hfill\vrule height6pt width6pt depth0pt}
\newcommand\ip[2]{\langle #1,#2\rangle}
\newcommand\cover[1][\Omega]{\tilde{#1}}
\newcommand\tfold[1][\Omega]{\tilde{#1}}
\newcommand\univ[1][\Omega]{\tilde{#1}_\infty}
\newcommand\q{\mathcal{Q}}
\newcommand\op[1][A]{H_{#1,V}^{}}
\newcommand\liftop[1][0]{\tilde H_{#1,V}^{}}
\newcommand\qf[1][A]{Q_{#1,V}^{}}
\newcommand\liftqf[1][0]{\tilde Q_{#1,V}^{}}
\newcommand\curl{{\rm \;curl\;}}
\newcommand\conj[1]{\overline{#1}}
\newcommand\nod[1][u]{\mathcal N(#1)}
\newcommand\schroedinger{\protect{Sch\-r\"o\-ding\-er}{} }
\newcommand\bohmaharonov{\protect{Bohm-Ahar\-o\-nov}{} }
\newcommand\slit{slit}

\newtheorem{theorem}{Theorem}[section]
\newtheorem{lemma}[theorem]{Lemma}
\newtheorem{corollary}[theorem]{Corollary}
\newtheorem{definition}[theorem]{Definition}
\newtheorem{proposition}[theorem]{Proposition}
\newtheorem{property}{Property}

\theorembodyfont{\rmfamily}
\newtheorem{example}[theorem]{Example}
\newtheorem{remark}[theorem]{Remark}
\newtheorem{remarks}[theorem]{Remarks}
\newtheorem{notation}[theorem]{Notation}

\makeatletter
\@addtoreset{equation}{section}
\makeatother

\newenvironment{proof}[1][]{\par\noindent\textit{Proof#1.} }{%
                               \vskip\baselineskip}%

\allowdisplaybreaks

\begin{document}
{\centering\LARGE
Nodal sets for the groundstate of the\\
\schroedinger operator with zero magnetic field\\
in a non simply connected domain.%
\footnote{Funded by the European Union TMR grant FMRX-CT 96-0001}\par
\vspace{1em}\large
B. Helffer$^1$\\
M. Hoffmann-Ostenhof$^2$\\
T. Hoffmann-Ostenhof$^{3,4}$\\
M. P. Owen$^4$
\par\vspace{1em}\normalsize
$^1$D\'epartement de Math\'ematiques, Universit\'e Paris-Sud\\ 
$^2$Institut f\"ur Mathematik, Universit\"at Wien\\
$^3$Institut f\"ur Theoretische Chemie, Universit\"at Wien\\
$^4$International Erwin \schroedinger Institute for Mathematical Physics
\par\vspace{1em}\today\par}

\begin{abstract}
  We investigate nodal sets of magnetic \schroedinger operators with zero
  magnetic field, acting on a non simply connected domain in $\r^2$. For the
  case of circulation $1/2$ of the magnetic vector potential around each
  hole in the region, we obtain a characterisation of the nodal set, and
  use this to obtain bounds on the multiplicity of the groundstate. For the
  case of one hole and a fixed electric potential, we show that the first
  eigenvalue takes its highest value for circulation $1/2$.
\end{abstract}

\section{Introduction and statement of results}

Let $\Omega\subset\r^2$ be a region with smooth ($C^\infty$) boundary, which
is homeomorphic to a disk with $k$ holes, and consider the magnetic
\schroedinger operator%
%
% The command
%
                        $\phantom{\mathbf{}}$%
%
% is inserted because of a strange bug in \LaTeX{}, where
% if the first occurence of either \mathcurl or \mathbf is not within the $$
% environment then \mathcurl will fail after the first occurence of \mathbf!
%
\begin{equation}
  \op:=(i\grad+A)^2+V
\end{equation}
acting on $L^2(\Omega)$ with Neumann boundary conditions. The potential $V$ is
assumed to be smooth, and we consider a smooth magnetic vector
potential $A$ which corresponds to a zero magnetic field. That is,
\begin{equation}\label{eqn:curl0}
  B:=\curl A=0
\end{equation}
in $\Omega$. Assumption~\eqref{eqn:curl0} implies that in any simply
connected, open subset of $\Omega$, there exists a gauge function $\phi$ such
that
\begin{equation}\label{eqn:Agrad}
  \grad\phi=A.
\end{equation}

We shall see that the operator $\op$ is unitarily equivalent to the
non-magnetic \schroedinger operator $\op[0]$ if and only if one can extend this
local gauge $e^{i\phi}$ to a globally defined function such that $\phi$ (which
might not be a singlevalued function) satisfies~\eqref{eqn:Agrad}. We shall see
that this can be done precisely when each of the circulations
\begin{equation}
  \Phi_i=\frac1{2\pi}\oint_{\sigma_i} A\cdot\d\x,
\end{equation}
of $A$ round the $i$-th hole ($i=1,\dots,k$) takes an integer value. Here
$\sigma_i$ is a closed
path\footnote{A piecewise smooth mapping $\gamma:[0,1]\rightarrow X$ is called
a path in $X$. The point $\gamma(0)$ is called the initial point and
$\gamma(1)$ is called the final point. The image $\Gamma=\gamma([0,1])$ of the
path is called a curve.} which parametrises the boundary $\Sigma_i$ of the
$i$-th hole and turns once in an anti-clockwise direction.

Furthermore, if the circulations $\Phi=(\Phi_1,\dots,\Phi_k)$ of two
distinct vector potentials $A$ and $A'$ are equal modulo
$\z^k$ then the corresponding operators $\op$ and $H_{A',V}$ are unitarily
equivalent under a gauge transformation.

\begin{theorem}\label{thm:properties}
  Let $\Omega\subset\r^2$ be a region with smooth boundary, which is
  homeomorphic to a disk with $k$ holes. For a given smooth potential $V$, the
  first
  eigenvalue $\lambda_1$ of the magnetic \schroedinger operator $\op$, where
  $A$ satisfies~\eqref{eqn:curl0}, depends only on the
  circulations $\Phi=(\Phi_1,\dots,\Phi_k)$ of $A$.
  The function $\lambda_1(\Phi)$ has the following properties (in which
  $l\in\z^k$ is arbitrary): 
  \begin{gather}
    \lambda_1(\Phi+l)=\lambda_1(\Phi)\label{eqn:first}\\
    \lambda_1(l/2+\Phi)=\lambda_1(l/2-\Phi)\label{eqn:second}\\
    \lambda_1(\Phi)>\lambda_1(0,\dots,0)\quad\text{for }\Phi\not\in\z^k.
    \label{eqn:third}
  \end{gather}
  For the case $k=1$, we have in addition to equation~\eqref{eqn:third} that
  \begin{equation}\label{eqn:fourth}
    \lambda_1(\Phi)<\lambda_1(1/2)
  \end{equation}
  for $\Phi\not\in1/2+\z$.
%and $\lambda_1(\Phi)=\lambda_1(1/2)$ if and only if $\Phi\in1/2+\z$.
\end{theorem}

Equations~\eqref{eqn:first}, \eqref{eqn:second} and
inequality~\eqref{eqn:third} are straightforward, and are proved in
Section~\ref{sec:elementary} (see also Remark~\ref{rem:gauge}). We choose
Neumann
boundary conditions on $\op$ in this article because we were motivated by
questions arising in the Ginzburg model of super-conductivity. Our
results are also valid for the case of Dirichlet boundary conditions (see
Remark~\ref{rem:mainthm}~\ref{item:neudir}). Dirichlet boundary conditions
are related to the \bohmaharonov effect for bounded states.
See~\cite{LaviOCar,Helf1,Helf2,Helf3}. Such models also arise in the
description of the Little-Parks experiment~\cite{LittPark}.

Inequality~\eqref{eqn:fourth} appears, to the best of our knowledge, for the
first time. Our proof of this result (see Section~\ref{sec:character}), uses a
connection between the maximality of the first eigenvalue for flux
$1/2$ and the structure of the nodal set of groundstates. The
nodal sets for the single hole case with flux $1/2$ were
recently investigated
by Berger and Rubinstein~\cite{BergRubi}. Part of our work is motivated by
their preprint.

Using semiclassical arguments as in~\cite{Helf1}, we can show that
in general the first eigenvalue is not necessarily maximised for circulation
$(1/2,\dots,1/2)$.

\begin{definition}
  The nodal set $\nod$ of an eigenfunction $u$ of a magnetic
  \schroedinger operator on a manifold $\Omega$ with smooth boundary is defined
  in $\overline\Omega$ by
  \begin{equation}
    \nod:=\conj{\{x\in\Omega:u(x)=0\}}.
  \end{equation}
\end{definition}

Some useful information on nodal sets of real valued eigenfunctions of
non-magnetic \schroedinger equations in two dimensions is given in
Proposition~\ref{thm:nodinfo}. In particular we see that such nodal sets
consist of the finite union of smoothly immersed circles and lines.
It is ``generically'' the case that the nodal set of every complex
eigenfunction of a magnetic
\schroedinger operator consists of isolated points of intersection of the
lines of zeros of the real and imaginary parts of the function.
See~\cite{ElliMataQi}.

The local properties of the nodal sets of eigenfunctions of the operator $\op$
are the same as the local properties of complex solutions of non-magnetic
\schroedinger
equations. More precisely, since we may find at every point a local gauge
$e^{i\phi}$ satisfying~\eqref{eqn:Agrad}, we may multiply any eigenfunction of
$\op$ by a local gauge so that the product solves a non-magnetic
\schroedinger equation. The nodal set is invariant under local gauge
transformations.

We shall see in what follows that although the local properties of nodal sets
of eigenfunctions of our magnetic \schroedinger operator are the same as the
properties of a non-magnetic \schroedinger operator, the global properties
differ in
the case where $\Phi=(1/2,\dots,1/2)$. In particular, in the non-magnetic
case we see that (since a real eigenfunction must change sign at the nodal set)
an even number of nodal lines (or perhaps no nodal lines) of an eigenfunction
emerges from each boundary component of the region. In
Theorem~\ref{thm:main} we show that for $\Phi=(1/2,\dots,1/2)$,
an odd number of nodal lines of the groundstate emerge from each component.

\begin{definition}
  We say that a (nodal) set $\mathcal N$ \textbf{\slit{s}} $\overline\Omega$
  if it is the union of a collection of piecewise smooth, immersed lines such
  that
  \begin{enumerate}\setlength\itemsep{0pt}
    \item each line starts and finishes at the boundary $\partial\Omega$ and
          leaves the boundary transversally;
    \item internal intersections between lines are transversal;
    \item\label{item:oddno} the complement $\Omega\setminus\mathcal N$ is
          connected;
    \item an odd number of nodal lines leaves each interior boundary component.
  \end{enumerate}
  We shall say that a collection of paths \slit{s} $\overline\Omega$ if the
  union of the images of the paths \slit{s} $\overline\Omega$.
\end{definition}

\singdiag{%
  \renewcommand\afigwidth{80}
  \aboundingbox{0}{0}{606}{458}}{%
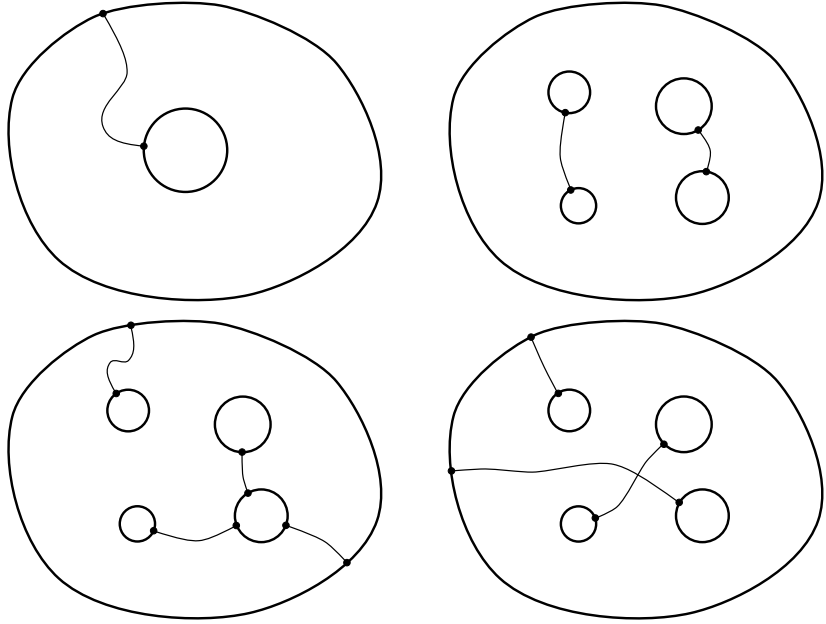
  \fig{cutting.eps}{\label{fig:cutting}}{}{Examples of some sets which \slit{} 
       $\overline\Omega$}}%
See Figure~\ref{fig:cutting} for some examples of regions which are \slit{}.
Note that part~\ref{item:oddno} of the above definition is the reason why a
nodal set which \slit{s} $\overline \Omega$ contains no immersed circles, and
also implies that each line of a \slit{ting} set links together a unique pair
$\{\Sigma_i,\Sigma_j\}$ of distinct (i.e.
$i\neq j$) boundary components. Note also that for the single hole case, a set
which \slit{s} $\overline\Omega$ consists of one line which joins the outer
boundary of $\Omega$ to the inner boundary.

In Corollary~\ref{thm:subsuper} we show that if a collection of paths \slit{} a
region then no sub- or supercollection of these paths can also \slit{} the
region. In Proposition~\ref{thm:boundlines} we show that the number $n$ of paths of
such a collection must satisfy $k/2\le n\le k$.

\begin{theorem}\label{thm:main}
  Let $\Omega$ be a region with smooth boundary, which is homeomorphic to a
  disk with $k$ holes. Let $V$ be a smooth potential and let $A$ be a smooth
  magnetic vector potential satisfying equation~\eqref{eqn:curl0}, such
  that the value of the circulations around each hole lie in $1/2+\z$ (that is
  $\Phi=(1/2,\dots,1/2)$, modulo $\z^k$). 
  \begin{enumerate}
    \item\label{item:cutting} If the first eigenvalue of $\op$ is simple then
    the nodal set of the corresponding eigenfunction \slit{s} $\overline\Omega
    $. Otherwise there exists an orthonormal
    basis $\{u_1,\dots,u_m\}$ of the groundstate ei\-gen\-spa\-ce such that
    the nodal set of any
    combination $\sum_{i=1}^ma_iu_i$, with $a_i\conj a_j\in\r$ for each $1\le
    i,j\le m$, \slit{s} $\overline\Omega$.
    \item\label{item:mult} The multiplicity $m$ of the first eigenvalue of         $\op$ satisfies
    \begin{equation}\label{eqn:mult}
      m\le\begin{cases}
        2,&k=1,2;\\
        k,&\text{$k$ odd, $k\ge3$};\\
        k-1,\quad&\text{$k$ even, $k\ge4$}.
      \end{cases}
    \end{equation}
    \item\label{item:emptynodal}
    For $k=1,2$ with groundstate multiplicity two, the nodal sets of
    two linearly independent groundstates do not intersect. It follows that the
    nodal set of a combination $a_1u_1+a_2u_2$ is empty whenever $a_1\conj a_2
    \not\in\r$.
  \end{enumerate}
\end{theorem}

Here we make some remarks connected to the above theorem.

\begin{remarks}\label{rem:mainthm}
\begin{enumerate}
\item The above bound on the multiplicity of the first eigenvalue is sharp in
the case of one hole (see Example~\ref{exa:sharp}), but it is not expected to
be sharp for many holes. It would be interesting to know an asymptotic result
about the growth of the maximum multiplicity with the number of holes.

\item We prove the bound by taking advantage of topological obstructions to
nodal sets caused by the holes. These obstructions prevent the existence of
high dimensional groundstate eigenspaces. Our type of
method was first discovered in~\cite{Chen} and has since been taken up and used
by others, e.g.~\cite{Nadi,HoffHoffNadi,HoffMichNadi}. See also~\cite{Coli} for
explicit constructions of examples with high multiplicity.

\item Our result bears
similarities to bounds on multiplicities of higher eigenvalues of non-magnetic
\schroedinger operators on surfaces with boundary. Some related literature on
this topic is given in~\cite{Coli,Nadi,HoffHoffNadi,HoffMichNadi}. 

\item It has been shown in~\cite{BessColbCour} that no upper bound on the
multiplicity exists when one adds a general magnetic field, even on the sphere.

\item For the cases $k\ge3$ we expect that there
could be intersection of nodal sets of two independent groundstates, and
correspondingly that the nodal set of a combination $a_1u_1+a_2u_2$ will not in
general be empty when $a_1\conj a_2\not\in\r$.

\item\label{item:neudir}
  If we assume that $\op$ has Dirichlet boundary conditions then
  Theorems~\ref{thm:properties} and~\ref{thm:main} hold with suitable changes
  to the proofs. More precisely, in Proposition~\ref{thm:nodinfo} the Taylor
  expansion~\eqref{eqn:cosexp} for a zero of order $l$ at a point $x\in\partial
  \Omega$ becomes
  \[ f(x)=ar^l\sin l\omega+O(r^{l+1}), \]
  and from Lemma~\ref{thm:bdzero} through to the proof of
  Theorem~\ref{thm:main}~\ref{item:mult}, all arguments which involve a
  function which has a zero of order $l=k$ (for example) should be replaced by
  the same argument involving a function with a zero of order $l=k+1$.
\end{enumerate}
\end{remarks}

\section{Some basic results}\label{sec:elementary}

The quadratic form corresponding to the operator $\op$ is
\begin{equation}\label{eqn:quaddef}
  \qf(u)=\int_\Omega\left(|(i\grad+A)u|^2+V|u|^2\right)\d^2x,
\end{equation}
with domain $\q^{\operatorname{Neu}}=W^{1,2}(\Omega)=H^1(\Omega)$. This choice
of quadratic form domain corresponds to Neumann boundary conditions for $\op$.
For the case of Dirichlet boundary conditions (see
Remark~\ref{rem:mainthm}~\ref{item:neudir}) the relevant quadratic form domain
is $\q^{\operatorname{Dir}}=W_0^{1,2}(\Omega)$.

\begin{remark}
  Neumann boundary conditions for a magnetic \schroedinger operator mean that
  functions in the domain of the operator satisfy
  \begin{equation}
    i\frac{\partial u}{\partial n}=-A\cdot n\;u
  \end{equation}
  on $\partial\Omega$, where $n$ is normal to $\partial\Omega$.

  One can always assume that the vector potential satisfies the additional
  properties
  \begin{equation}\label{eqn:divfree}
    \grad\cdot A=0 \text{ in }\Omega,\qquad\qquad A\cdot n=0\text{ on }\partial
    \Omega.
  \end{equation}
  The reason is as follows: There is a solution $\phi$ (unique up to a
  constant) to the oblique derivative problem
  \begin{equation}\Delta\phi=-\grad\cdot A\text{ in }\Omega,\qquad\qquad\grad
  \phi\cdot n=-A\cdot n\text{ on }\partial\Omega.
  \end{equation}
  See~\cite[Theorem 6.31 and the following remark]{GilbTrud}.
  Setting $A'=A+\grad\phi$, the operator $\op[A']$ is
  unitarily equivalent to $\op$ under the gauge transformation $e^{i\phi}$,
  and $A'$ satisfies the properties~\eqref{eqn:divfree}.
\end{remark}

\begin{proof}[ of equation~\eqref{eqn:first}]
  Let $A$ and $A'$ be magnetic vector potentials with circulations that differ
  by an element of $\z^k$. For any closed path $\sigma$,
  \[ \frac1{2\pi}\oint_\sigma (A'-A)\cdot\d\x\in\z, \]
  and hence there exists a smooth, multivalued function $\phi$ such that
  $e^{i\phi}$ is univalued and $\grad\phi=A'-A$. For $u\in H^1(\Omega)$ we have
  \[ (i\grad+A')e^{i\phi}u=e^{i\phi}(i\grad+A)u, \]
  and therefore the operators $\op$ and $\op[A']$ are unitarily equivalent.\qed
\end{proof}

\begin{remark}\label{rem:gauge}
  For any magnetic vector potential $A$ satisfying~\eqref{eqn:curl0} there
  exists a gauge function $\phi$ such that
  \[ A(x,y)-\sum_{i=1}^k\frac{\Phi_i}{2\pi r_i^2}\binom{-y+y_i}{x-x_i}=(\grad
     \phi)(x,y), \]
  where $(x_i,y_i)$ is a fixed point in the $i$-th hole, $r_i^2=(x-x_
  i)^2+(y-y_i)^2$ and $\Phi_i$ is the circulation of $A$ round the
  $i$-th hole. Defining
  \[ A'(x,y)=\sum_{i=1}^k\frac{\Phi_i}{2\pi r_i^2}\binom{-y+y_i}{x-x_i}. \]
  we see, for a fixed $V$, that
  \[ \op[A']=e^{-i\phi}\op e^{i\phi} \]
  and thus $\op$ is unitarily equivalent to $\op[A']$. This means
  that the magnetic vector potential is determined up to a gauge transformation
  by its circulations $\Phi$, and verifies that the spectrum of $\op$ is
  determined by $\Phi$.
\end{remark}

\begin{proof}[ of equation~\eqref{eqn:second}]
  Let $A$ be a magnetic vector potential with circulation $\Phi$, and let $u$
  be a groundstate of $\op$. It is easy to show that $\overline u$ is a
  groundstate of $\op[-A]$ with the same eigenvalue, and hence
  \begin{equation}\label{eqn:erm}
    \lambda_1(-\Phi)=\lambda_1(\Phi).
  \end{equation}
  We obtain equation~\eqref{eqn:second} by combining~\eqref{eqn:erm}
  and~\eqref{eqn:first} as follows:
  \begin{equation*}
    \lambda_1(l/2+\Phi)=\lambda_1(-l/2-\Phi)=\lambda_1(l/2-\Phi).\tag*{\qed}
  \end{equation*}
\end{proof}

\begin{proof}[ of inequality~\eqref{eqn:third}]
  Suppose for a contradiction that $\Phi\not\in\z^k$ and that $\lambda_1(\Phi)
  \le\lambda_1(0)$, where $\Phi$ is the circulation vector of some magnetic
  vector potential $A$. Let $u_0$ denote the unique normalised
  positive groundstate of the operator $\op[0]$ and let $u_A$ be a normalised
  groundstate of the operator $\op$. Using the diamagnetic
  inequality~\cite{Simo} we have
  \begin{equation}\label{eqn:kato}
    \qf[0](|u_A|)\le\qf(u_A)=\lambda_1(\Phi)\le\lambda_1(0)=\qf[0](u_0),
  \end{equation}
  and thus $|u_A|=u_0$. It follows that $u_A=e^{i\phi}u_0$ for some smooth,
  real valued, multivalued function $\phi$, and hence
  \begin{align*}
    \int_\Omega|A-\grad\phi|^2|u_0|^2\d^2x&=\int_\Omega|(i\grad+A-\grad\phi)u_0
    |^2\d^2x-\int_\Omega|\grad u_0|^2\d^2x\\
    &=\int_\Omega|(i\grad+A)u_A|^2\d^2x-\int|\grad u_0|^2\d^2x\\
    &=\qf(u_A)-\qf[0](u_0)\\
    &=0,
  \end{align*}
  and therefore $A=\grad\phi$ in $\Omega$. Thus for each $i=1,\dots,k$ we have
  \[ \Phi_i=\frac1{2\pi}\oint_{\sigma_i}A\cdot\d\x=\frac1{2\pi}\oint_{\sigma_i}
     \d\phi\in\z, \]
  where $\sigma_i$ is a closed path which parametrises the boundary $\Sigma_i$
  of the $i$-th hole and turns once in an anticlockwise direction.
  This contradicts our assumption that $\Phi\not\in\z^k$.\qed
\end{proof}

The proof of inequality~\eqref{eqn:third} is an alternative to the proofs given
in~\cite{LaviOCar} and~\cite{Helf1}. It has the advantage of being simpler and being
independent of whether the boundary conditions are Neumann or Dirichlet.

We leave the proof of inequality~\eqref{eqn:fourth} until
Section~\ref{sec:character} because it depends on
Theorem~\ref{thm:main}~\ref{item:cutting}.

\section{A twofold Riemannian covering manifold}

In this section we consider the case where the circulations of the magnetic
vector potential $A$ satisfy
\begin{equation}\label{eqn:halfcirc}
  \Phi_i\in1/2+\z
\end{equation}
for each $1\le i\le k$.
The proofs of our results use a twofold Riemannian covering
manifold $\tfold$ of the domain $\Omega$ (see Remark~\ref{rem:K} however). For
the case of more than one hole,
there exists more than one twofold Riemannian covering manifold of $\Omega$.
We shall take a particular choice of covering manifold on which the circulation
of the lifted magnetic ($1$-form) potential $\tilde A$ along any closed curve
is an integer. Before the precise definition, we introduce some basic notation.
For further details see for example~\cite{Kosn} or~\cite{GallHuliLafo}.

\begin{notation}
  Let $\tfold$ be a covering manifold of $\Omega$, and let $\Pi$ be the
  associated covering map. We denote the lifts of various quantities as
  follows:

  For a set $\mathcal N$ define $\cover[\mathcal N]=\{x\in\tfold:\Pi(x)\in
  \mathcal N\}$. For a function $f:\Omega\rightarrow\c$, define $\cover[f]:
  \tfold\rightarrow\c$ by $\cover[f]=f\circ\Pi$. For a path $\sigma:[0,1]
  \rightarrow\Omega$ and a point $x\in\tfold$ such that $\Pi(x)=\sigma(0)$
  let $\cover[\sigma]:[0,1]\rightarrow\tfold$ denote the unique lifted path
  such that $\cover[\sigma](0)=x$ and $\Pi\circ\cover[\sigma]=\sigma$.

  We endow the covering manifold with the metric obtained by lifting the
  flat Euclidean metric of $\Omega$ to $\tfold$. This is the unique metric
  which makes $\Pi$ a local isometry, and therefore a Riemannian covering
  map.
  Let $\cover[\Delta]=\operatorname{div}\operatorname{grad}$ denote the
  Laplace-Beltrami operator on $L^2(\tfold)$ induced by the lifted metric
  on $\Omega$, and let $\tilde A$ be the $1$-form on $\tfold$ obtained
  by lifting the $1$-form associated with the smooth vector potential $A$
  defined on $\Omega$.
\end{notation}

Let $\univ$ be the universal covering manifold of $\Omega$ and let $\Pi_\infty$
be the associated covering map. The universal covering of any manifold is
simply connected.

Note that due to~\eqref{eqn:halfcirc} if two points $x_\infty,y_\infty\in\univ$
satisfy $\Pi_\infty(x_\infty)=\Pi_\infty(y_\infty)$ then for any path $\sigma$
joining $x_\infty$ to $y_\infty$, the integral
\begin{equation}\label{eqn:equivclass}
  \frac1{2\pi}\oint_{\Pi_\infty\circ\sigma} A\cdot\d\x
\end{equation}
lies either in $1/2+\z$ or in $\z$. The value of~\eqref{eqn:equivclass} is
independent of the path $\sigma$ because $\curl A=0$ and because the universal
covering manifold is simply connected. We therefore construct the twofold
covering manifold (as a quotient of the universal covering manifold) as
follows:

\begin{definition}\label{def:tfold}
  \begin{enumerate}\item
  We define the twofold covering manifold $\tfold$ by identifying points $x_\infty$,
  $y_\infty$ in $\univ$ according to the equivalence relation $x_\infty\sim y_\infty$ if and only if
  \begin{equation}\label{eqn:above}
    \Pi_\infty(x_\infty)=\Pi_\infty(y_\infty)
  \end{equation}
  and for each path $\sigma$ in $\univ$ joining $x_\infty$ to $y_\infty$ we
  have
  \begin{equation}\label{eqn:intflux}
    \frac1{2\pi}\int_{\Pi_\infty\circ\sigma} A\cdot\d\x\in\z.
  \end{equation}
  The covering map $\Pi:\tfold\rightarrow\Omega$ is defined by $\Pi(x)=\Pi_
  \infty(x_\infty)$, where $x=[x_\infty]$ is the equivalence class (under
  $\sim$) containing $x_\infty$.
\singdiag{%
  \renewcommand\afigwidth{50}%
  \aboundingbox{0}{0}{590}{493}}{%
  \fig{altnewcover.eps}{\label{fig:twocover}}{%
  \put(190,50){$\Pi$}%
  \put(-20,20){$\Omega$}%
  \put(-20,115){$\tfold$}}%
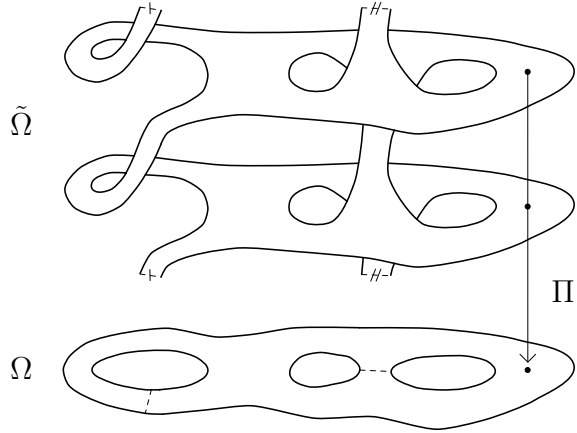
  {Realization of a twofold covering manifold}}%
\item On our twofold covering manifold we define the symmetry map
$\symm:\tfold\rightarrow\tfold$ by setting $\symm x$ to be the other point in
$\tfold$ which lies above $\Pi(x)\in\Omega$. Note that $\Pi^{-1}(\Pi(x))=\{x,
\symm x\}$.
\item
  We say that a function $f:\tfold\rightarrow\c$ is symmetric if $f(\symm x)=f(
  x)$ for all $x\in\tfold$, and antisymmetric if $f(\symm x)=-f(x)$ for all
  $x\in\tfold$.
\end{enumerate}
\end{definition}

Note that the identity map and
$\symm$ form a group $\mathcurl G=\{I,\symm\}$, with the composition $\symm^2=I
$, which acts freely on $\tfold$. The quotient of $\tfold$ by $\mathcurl G$ is
the original manifold $\Omega$.
The lift $\cover[f]$ of a function $f$ on $\Omega$ is symmetric.

Using equation~\eqref{eqn:intflux} we have
\begin{equation}
  \frac1{2\pi}\oint_\sigma\tilde A\cdot\d\tilde{\x}=\frac1{2\pi}\oint_{\Pi\circ
  \sigma}A\cdot\d\x\in\z,
\end{equation}
for any closed path $\sigma$ in $\tfold$. Hence there exists a
smooth, multivalued function $\theta$ on $\tfold$ such that $\exp{i\theta}$ is
univalued and
\begin{equation}\label{eqn:gradeqA}
  \operatorname{grad}\theta=\tilde A.
\end{equation}

\begin{lemma}\label{thm:L}
  The operator $\mathcurl L:L^2(\Omega)\rightarrow L^2(\tfold)$ defined by
  \begin{equation}\label{eqn:defL}
    \mathcurl Lu=\frac1{\sqrt2}e^{i\theta}\tilde u
  \end{equation}
  is a isometry onto the antisymmetric functions in $L^2(\tfold)$, and maps
  eigenfunctions of $\op$ onto antisymmetric eigenfunctions of the
  \schroedinger operator
  \begin{equation}\label{eqn:plliftop}
    \liftop=-\tilde\Delta+\tilde V
  \end{equation}
  acting on $L^2(\tfold)$ with Neumann boundary conditions. 
\end{lemma}

\begin{proof}
  We shall first show that the function $e^{i\theta}$ is antisymmetric (under
  $\symm$). For any point $x\in\tfold$, let
  $\sigma:[0,1]\rightarrow\tfold$ be a path which joins $x$ to $\symm x$.
  Using the terminology of Definition~\ref{def:tfold} we have
  $\Pi(x)=\Pi(\symm x)$ but $x\not\sim Gx$, and hence
  \[ \frac1{2\pi}\oint_{\Pi\circ\sigma}A\cdot\d\x=l+1/2 \]
  for some $l\in\z$. Keeping in mind that $\theta$ is multivalued, we get
  \begin{equation*}
    \theta(\symm x)-\theta(x)=\int_\sigma\d\theta
    =\int_\sigma\tilde A\cdot\d\tilde\x 
    =\oint_{\Pi\circ\sigma}A\cdot\d\x 
    =(2l+1)\pi.
  \end{equation*}
  Hence $\exp[i\theta(\symm x)]=-\exp[i\theta(x)]$ as claimed.

  The action of $\mathcurl L$ upon a function $u\in L^2(\Omega)$ consists of
  two steps. The first step is to lift $u$ to the symmetric function
  $\tilde u$. This is a bijection onto the space of symmetric functions of
  $L^2(\tfold)$. The second step is to multiply $\tilde u$ by the antisymmetric
  function $e^{i\theta}$. This step is a bijection from the space of symmetric
  functions onto the space of antisymmetric functions in $L^2(\tfold)$. To see
  that $\mathcurl L$ is an isometry onto its range, we take two functions
  $u,v\in L^2(\Omega)$ and note that
  \[ \ip{\mathcurl Lu}{\mathcurl Lv}_{L^2(\tfold)}=\frac12\int_{\tfold}e^{i\phi
     }\tilde u.e^{-i\phi}\overline{\tilde v}\d\tilde x=\int_\Omega u\overline
     v\d x=\ip uv_{L^2(\Omega)}. \]

  For every eigenfunction $u$ of $\op$, the lift $\tilde u$ is an eigenfunction
  of the lifted magnetic \schroedinger operator
  \begin{equation}
    \liftop[A]=(i\operatorname{div}+\tilde A)(i\operatorname{grad}+
    \tilde A)+\tilde V
  \end{equation}
  on $\tfold$ where $\tilde V$ and $\tilde A$ are the lifts of $V$ and 
  $A$ respectively. We now multiply by the gauge $e^{i\theta}$. Using
  equation~\eqref{eqn:gradeqA}, the function $e^{i\theta}\tilde u$ is an
  eigenfunction of the non-magnetic \schroedinger operator $\liftop$.\qed
\end{proof}

The spectrum of $\op$ consists of the eigenvalues
corresponding to the antisymmetric eigenfunctions of $\liftop$. It turns out to
be useful (see Lemma~\ref{thm:nodal}) to single out the case where a function
$u$ has the following property:

\setcounter{property}{15}
\begin{property}\label{pro:reallift}
  The function $u$ is a groundstate of the operator $\op$, and
  the corresponding eigenfunction $\mathcurl Lu$ of $\liftop$ has a constant
  phase. In other words, there exists a constant $\alpha\in\c\setminus\{0\}$
  such that $\mathcurl L(\alpha u)$ is a real valued function.
\end{property}

Due to the symmetry of $\tfold$, the groundstate of the operator $\liftop$ is
symmetric. In contrast, if $u$ has Property~\ref{pro:reallift} then $\mathcurl
L(\alpha u)$ is an antisymmetric eigenfunction (and therefore an excited state)
of $\liftop$. Consequently both $\mathcurl L(\alpha u)$ and $u$ have a nonempty
nodal set.

\begin{remark}\label{rem:K}
  It is not necessary to use the covering manifold to describe
  Property~\ref{pro:reallift}. An alternative is to formulate the property in
  terms of an antilinear operator $K$. We define the operator below.

  Since $\Phi_i\in1/2+\z$ for each $i=1,\dots,k$, we see that
  \[ \frac1{2\pi}\oint_\sigma 2A\cdot\d\x\in\z \]
  for all closed paths $\sigma$ in $\Omega$. It follows that there exists a
  smooth, multivalued function $\psi$ such that $e^{i\psi}$ is univalued and
  $\grad\psi=2A$. The multivalued function $\theta$ given in
  equation~\eqref{eqn:gradeqA} is related to $\psi$ by the
  formula
  \[ \psi\circ\Pi=2\theta+c \]
  for some constant $c$.
  We define $K$ by the formula
  \begin{equation}
    K=e^{-i\psi}\Gamma,
  \end{equation}
  where $\Gamma$ is the operator $\Gamma u=\overline u$. Then $K^2=
  \operatorname{Id}$ and $K$ commutes with $\op$.
  It turns out that a function $u\in L^2(\Omega)$ has
  Property~\ref{pro:reallift} if and only if it is an eigenfunction of both
  $\op$ and $K$.

  One could in fact completely dispense with the covering manifold, but at the
  expense of a clear geometrical picture in the following sections.
\end{remark}

\section{Characterisation of the nodal set}\label{sec:character}

We first collect some well known facts about eigenfunctions
of non-magnetic \schroedinger operators acting on two dimensional Riemannian
manifolds:
\begin{proposition}[Non-magnetic \schroedinger operators]\label{thm:nodinfo}
  Let $f$ be a real valued eigenfunction of a non-magnetic \schroedinger
  operator with smooth potential and Neumann boundary conditions, on a two
  dimensional locally flat Riemannian manifold $\Omega$ with smooth boundary.
  Then $f\in C^\infty(\overline\Omega)$. Furthermore, $f$ has the following
  properties:
  \begin{enumerate}
    \item\label{item:taylor} If $f$ has a zero of order $l$ at a point $x_0\in
          \overline\Omega$ then the Taylor expansion of $f$ is
          \begin{equation}\label{eqn:harmpoly}
            f(x)=p_l(x-x_0)+O(|x-x_0|^{l+1}),
          \end{equation}
          where $p_l$ is a real valued, non-zero, harmonic, homogeneous
          polynomial of degree $l$.

          Moreover if $x_0\in\partial\Omega$, the Neumann boundary
          conditions imply that
          \begin{equation}\label{eqn:cosexp}
            f(x)=ar^l\cos l\omega+O(r^{l+1})
          \end{equation}
          for some non-zero $a\in\r$, where $(r,\omega)$ are polar coordinates
          of $x$ around $x_0$. The angle $\omega$ is chosen so that the tangent
          to the boundary at $x_0$ is given by the equation $\sin\omega=0$.
    \item \label{item:nod}The nodal set $\nod[f]$ is the union of finitely many, smoothly
          immersed circles in $\Omega$, and smoothly immersed lines which
          connect points of
          $\partial\Omega$. Each of these immersions is called a
          \textit{nodal line}. Note that self-intersections are allowed.
          The connected components of $\Omega\setminus\nod[f]$ are called
          \textit{nodal domains}.
    \item \label{item:order} If $f$ has a zero of order $l$ at a point $x_0\in
          \Omega$ then exactly $l$ segments of nodal lines pass through $x_0$.
          The tangents to the nodal lines at $x_0$ dissect the full circle into
          $2l$ equal angles.

          If $f$ has a zero of order $l$ at a point $x\in
          \partial\Omega$ then exactly $l$ segments of nodal lines meet the
          boundary at $x_0$. The tangents to the nodal lines at $x_0$ are given
          by the equation $\cos l\omega=0$, where $\omega$ is chosen as
          in~\eqref{eqn:cosexp}.
  \end{enumerate}
\end{proposition}

\begin{proof}
  The proof that $f\in C^\infty(\overline\Omega)$ can be found 
  in~\cite[Theorem 20.4]{Wlok}.

  The proof of part~\ref{item:taylor} is trivial because $V$ and $f$ are
  smooth functions so the Taylor expansion (with remainder) exists. The
  properties of the
  first term of the expansion follow by substituting the Taylor expansion into
  the groundstate eigenvalue equation. 

  See~\cite{Bers},~\cite{Chen} for proofs of the other parts.\qed
\end{proof}

Proposition~\ref{thm:nodinfo} can be generalised to include eigenfunctions
of magnetic \schroedinger operators with a smooth magnetic vector potential
$A$. The eigenfunctions still lie in $C^\infty(\overline\Omega)$ and the
expansions~\eqref{eqn:harmpoly} and~\eqref{eqn:cosexp} hold, except that the
polynomial $p_l$ and the constant $a$ are allowed to be complex. However 
statements~\ref{item:nod} and~\ref{item:order} about the nodal set do not carry over.

\begin{theorem}\label{thm:equiv}
  Let $\mathcal N\subset\Omega$ be the union of finitely many smoothly immersed
  circles
  and smoothly immersed lines which connect points of $\partial\Omega$. The
  following statements are equivalent:
  \begin{enumerate}
    \item $\Omega\setminus\mathcal N$ is connected (therefore $\mathcal N$
          contains no smoothly immersed circles), and an odd number of lines
          emanate from each hole.
    \item In the twofold covering manifold, the open set $\tfold\setminus\tilde
          \mathcal N$ decomposes into two
          open path connected subsets $D_1,D_2$ such that $D_2=\symm D_1$ and
          $\partial D_1\cap\tfold=\partial D_2\cap\tfold=\tilde\mathcal N$.
  \end{enumerate}
\end{theorem}

\begin{proof}
  (i)$\Rightarrow$(ii) Let $D_1$ be a connected component of $\tfold\setminus
  \tilde\mathcal N$. Suppose for a contradiction that this is the only
  component.
  Due to the symmetry of the manifold, $\symm D_1=D_1$, and thus for any point
  $x\in D_1$ there exists a path $\sigma$ lying in $D_1$ (i.e. not intersecting
  $\tilde\mathcal N$), which joins $x$ and $Gx$. Using the terminology of
  Definition~\ref{def:tfold} we have $\Pi(x)=\Pi(\symm x)$ but $x\not\sim Gx$,
  and hence
  \[ \frac1{2\pi}\oint_{\Pi\circ\sigma}A\cdot\d\x\in1/2+\z. \]
  The closed path $\Pi\circ\sigma$ must therefore circulate an odd number of
  holes. Since an odd number of lines of $\mathcal N$ emanate from each hole,
  the path $\Pi\circ\sigma$ must intersect with one of them. This contradicts
  the fact that $\sigma$ does not intersect $\tilde\mathcal N$.

  Since $\Omega\setminus\mathcal N$ is connected there can only be two
  connected components $D_1,D_2$ of $\tfold\setminus\tilde\mathcal N$. As
  above, we see that $\symm D_1\neq D_1$, and therefore $D_2=\symm D_1$.

  Suppose now for a contradiction that $\partial D_1\cap\tfold\neq\tilde
  \mathcal N$. Then there
  exists a point $x\in\partial D_1\cap\tfold$ such that $x\not\in\partial D_2
  \cap\tfold$. The set $D_1$ borders with itself at $x$, and since $D_1$ is
  path connected there exists a closed path $\sigma$ such that $\sigma(0)=
  \sigma(1)=x$, which intersects $\tilde\mathcal N$ transversally at $x$ and
  which does not intersect $\tilde\mathcal N$ anywhere else. Since $\sigma$ is
  closed,
  \[ \oint_{\Pi\circ\sigma}A\cdot\d\x\in\z, \]
  and therefore $\Pi\circ\sigma$ circulates an even number of holes. Since an
  odd number of lines emanate from each hole, $\Pi\circ\sigma$ intersects
  $\mathcal N$ an even number of times. This contradicts the fact that $\sigma$
  intersects $\tilde\mathcal N$ only once.

  (ii)$\Rightarrow$(i) Since $D_2=\symm D_1$ we see that $\Pi D_1=\Pi D_2
  $, and hence $\Omega\setminus\mathcal N=\Pi(D_1\cup D_2)=\Pi D_1\cup\Pi D_2=
  \Pi
  D_1$. Since $\Pi$ is continuous, $\Omega\setminus\mathcal N$ is connected.
  Let $\sigma\subset\Omega$ be a closed path which circulates the $i$-th hole.
  Due to the construction of $\tfold$, $\sigma$ may be lifted to a path
  $\tilde\sigma$ in $\tfold$ which begins at a point $x\in D_1$ and ends at
  $Gx\in D_2$. Since $D_1$ and $D_2$ coborder, the path $\tilde\sigma$ crosses
  $\tilde\mathcal N$ an odd number of times and therefore $\sigma$ crosses
  $\mathcal N$
  an odd number of times. By choosing $\sigma\subset\Omega$ sufficiently close
  to $\sigma_i$ we see that an odd number of segments of lines leave the $i$-th
  boundary component. Since $\Omega\setminus\mathcal N$ is connected, each of
  these line endings belongs to a distinct line, and hence an odd number of
  lines leaves each boundary component.\qed
\end{proof}

\begin{corollary}\label{thm:subsuper}
  Suppose that a collection of paths \slit{s} a region. Then no subcollection
  of these paths can \slit{} the region. Also, no supercollection of these
  paths (i.e. a collection of paths which contain the original collection) can
  \slit{} the region. 
\end{corollary}

\begin{proof}
  Suppose that the union $\mathcal N$ of a collection
  of lines $\{\Gamma_1,\dots,\Gamma_n\}$ \slit{s} $\overline\Omega$. Using
  Theorem~\ref{thm:equiv}, we see that in the twofold covering manifold the
  open set $\tfold\setminus\tilde\mathcal N$ decomposes into two cobordering,
  open, path connected subsets $D_1,D_2$. Let $S$ be the union of a strict
  subcollection of the lines. The non-empty set $\tilde{\mathcal N}\setminus
  \tilde S$ connects together the two regions $D_1$ and $D_2$ and thus
  $\tfold\setminus\tilde S=D_1\cup D_2\cup(
  \tilde{\mathcal N}\setminus\tilde S)$ is connected. Using
  Theorem~\ref{thm:equiv} in the reverse direction, we see that $S$ does not
  \slit{} $\overline\Omega$.

  It follows easily that no supercollection of $\mathcal N$ can \slit{}
  because then $\mathcal N$ would be a strict subset of $S$ which
  \slit{s} $\overline\Omega$, and this is not possible by the above
  paragraph.\qed
\end{proof}

\begin{lemma}\label{thm:nodal}
  If a groundstate $u$ of $\op$ has Property~\ref{pro:reallift} then the
  nodal set of $u$ \slit{s} $\overline\Omega$.
\end{lemma}

\begin{proof} 
  By multiplying the function $u$ by a non-zero complex constant we may assume
  that the eigenfunction $\mathcurl Lu$ of $\liftop$ is real valued.
  Since $\mathcurl Lu$ is an
  antisymmetric function on the covering manifold $\tfold$, the nodal domains
  $D_1,\dots,D_l$ of $\mathcurl Lu$ have the property that for
  each $i=1,\dots,l$, we have $\symm D_i=D_j$ for some $j\neq i$. Suppose for a
  contradiction that
  $l>2$. Then there exist two cobordering domains $D_1,D_2$ such that $\symm D_
  1\neq D_2$. Define $D=\operatorname{Interior}(\overline{D_1\cup D_2})$, so
  that $D$ is the union of $D_1$, $D_2$ and the border between them.
  Let ${\tilde Q}^D_{0,V}$ denote the quadratic form corresponding to the
  \schroedinger operator
  \[ {\tilde H}_{0,V}^D=-\tilde\Delta+\tilde V \]
  on $D$ with Dirichlet boundary conditions on $\tilde S=\partial D
  \cap\tfold$ and Neumann boundary condition $\partial D\cap\partial\tfold$,
  and let $g$ denote the corresponding positive groundstate. Since the boundary
  of $D$ is piecewise smooth, the restriction $\mathcurl Lu|_D$ lies in the
  quadratic form domain of ${\tilde Q}^D_{0,V}$. Define the antisymmetric
  function $h$ on $\tfold$ by
  \[ h(y)=\begin{cases}g(y),&y\in D,\\-g(\symm y),&y
     \in\symm D,\\0,&\text{otherwise}.\end{cases} \]
  Let $\liftqf$ denote the quadratic form of the operator $\liftop$, which we
  define in equation~\eqref{eqn:plliftop}.
  Since $\mathcurl Lu$ is an antisymmetric eigenfunction which corresponds to
  a groundstate of $\op$, it has the least energy of all antisymmetric
  functions, and therefore
  \begin{equation}\label{eqn:ineq1}
    \frac{\liftqf(\mathcurl Lu)}{\|\mathcurl Lu\|_{L^2(\tfold)}^2}\le\frac{
    \liftqf(h)}{\|h\|_{L^2(\tfold)}^2}=\frac{{\tilde Q}^D_{0,V}(g)}{\|g\|_{L^2(
    D)}^2}\le\frac{{\tilde Q}^D_{0,V}(\mathcurl Lu|_D)}{\|\mathcurl Lu|_D
    \|_{L^2(D)}^2}=\frac{\liftqf(\mathcurl Lu)}{\|\mathcurl Lu\|_{L^2(\tfold)}^
    2}.
  \end{equation}
  We have in fact equality in~\eqref{eqn:ineq1}, and therefore, by uniqueness
  of the groundstate, we have that $\mathcurl Lu|_D=\lambda g$ for some
  $\lambda\neq0$. This contradicts the fact that $\mathcurl Lu|_D$ is zero on
  $\partial D_1\cap D$. Hence $l=2$ and the result follows from the equivalence
  proved in Theorem~\ref{thm:equiv}.\qed
\end{proof}

\begin{lemma}\label{thm:bdzero}
  If a groundstate $u$ of $\op$ has a zero of order $l$ at a point
  $x\in\partial\Omega$ then $l\le k$. Moreover, if $k$ is even and $x$ lies on
  an interior boundary component ($\Sigma_1$, say) then $l\le k-1$.
\end{lemma}

\begin{proof}
  Assume first that $u$ has Property~\ref{pro:reallift}, and suppose for a
  contradiction that $l\ge k+1$. Let $\Sigma_i$ denote the boundary component
  on which $x$ lies, where $i\in\{0,1,\dots,k\}$. At least
  $k+1$ distinct nodal lines emerge from $\Sigma_i$. Since there are only $k$
  boundary components distinct from $\Sigma_i$ there must exist two nodal lines
  which both start at $\Sigma_i$ and finish at $\Sigma_j$ for some $j\neq i$.
  In both cases, such a nodal set would split $\overline\Omega$ into more than
  one nodal domain, thus contradicting the assumption that $\nod$ \slit{s}
  $\overline\Omega$. Hence $l\le k$.

  If $u$ does not have Property~\ref{pro:reallift} then we can obtain a
  contradiction using the same methods above on the function $\mathcurl L^{-1}[
  \operatorname{Re}[\mathcurl Lu]]$. This function is a groundstate of $\op$,
  has a zero of order at least $l$ at $x$, and does have
  Property~\ref{pro:reallift}.

  Suppose that $k$ is even, that $x\in\Sigma_i$ (with $i\in\{1,\dots,k\}$) and
  that $l=k$. Since $\nod$ \slit{s} $\overline\Omega$ there must be an odd
  number of nodal lines leaving $\Sigma_i$. Therefore at least $k+1$ nodal
  lines leave $\Sigma_i$, and we obtain a contradiction as before.\qed
\end{proof}

\begin{proof}[ of Theorem~\ref{thm:main}~\ref{item:cutting}]
  Let $U$ denote the groundstate eigenspace of $\op$. For all $u\in U$ we have
  $\operatorname{Re}[\mathcurl Lu],\operatorname{Im}[\mathcurl Lu]\in\mathcurl
  LU$ are eigenfunctions of
  $\liftop$, if they are not identically zero. It follows that we may find an
  orthonormal basis $\{f_1,\dots,f_m
  \}$ of real valued functions for $\mathcurl LU$.  Since $\mathcurl L$ is an
  isometry, the functions $\{u_1,\dots,u_m\}$ defined by $u_i=\mathcurl L^{-1}f
  _i$ are an orthonormal basis of $U$.
  
  Now let $u=\sum_{i=1}^m\alpha_iu_i$, where $\alpha_i\overline\alpha_j\in\r$
  for each $1\le i,j\le m$. Take some $\alpha_j\neq0$. Then
  \[ \mathcurl L(\overline{\alpha_j}u)=\sum_{i=1}^m\alpha_i\overline{\alpha_j}
     f_i \]
  is a real valued function, and so $u$ has Property~\ref{pro:reallift}.
  The result now follows from Lemma~\ref{thm:nodal}.\qed
\end{proof}

\begin{lemma}\label{thm:zero}
  Suppose that the groundstate eigenspace $U$ of $\op$ is $m$ dimensional.
  \begin{enumerate}
    \item\label{item:mdim} For each point $x\in\partial\Omega$ there exists a
    function $u_x\in U$ which has Property~\ref{pro:reallift} and which has a
    zero of order at least $m-1$ at $x$.
    \item\label{item:unique} If $m=k+1$ then for each point $x$ lying on the
    outer boundary $\Sigma_0$ of $\overline\Omega$ there exists a unique $u_x
    \in U$ (up to multiplication by a complex constant) which has a zero of
    order $k$ at $x$. The function $u_x$ has Property~\ref{pro:reallift}. The
    nodal set of $u_x$ consists of $k$ lines which emanate from $x$ (which is
    the only point of intersection of lines), and which end at each of
    the $k$ distinct interior boundary components of $\Omega$. Each nodal line
    depends smoothly on $x$.
    \item\label{item:uni} If $k$ is even and $m=k$ then for each point $x$ 
    lying on an interior component of the boundary of $\overline\Omega$
    there exists a unique $u_x\in U$ (up to multiplication by a complex
    constant) which has a zero of order $k-1$ at $x$. The function $u_x$ has
    Property~\ref{pro:reallift}.
  \end{enumerate}
\doubdiag{%
    \renewcommand\afigwidth{40}%
    \renewcommand\bfigwidth{40}%
    \aboundingbox{0}{0}{298}{237}
    \bboundingbox{0}{0}{282}{225}}%
    {\fig{klines.eps}{\label{fig:klines}}{%
    \put(76,-8){$x$}}{}}%
    {\fig{kevelines.eps}{\label{fig:kevelines}}{%
    \put(62,105){$x$}}{}}%
\end{lemma}

\begin{proof}
  \ref{item:mdim} We shall first prove by induction the following
  statement: If $U_m$ is an $m$ dimensional vector space of groundstates of
  $\op$ then for each point $x\in\partial\Omega$ there exists a function $f\in
  U_m$ which has a zero of order at least $m-1$ at $x$.

  The first step of the induction, for $m=1$, is trivial. Assume now that the
  above statement is true for some general $m$.
  Suppose that $U_{m+1}$ is an $m+1$ dimensional vector space of eigenfunctions
  of $\op$. Let $U_m$ be any $m$ dimensional subspace of $U_{m+1}$. Then
  there exists a function $f_1\in U_m$ which has a zero of order at least $m-1
  $ at $x$. We can assume that the order of the zero is exactly $m-1$,
  otherwise we have found a function with a zero of order at least $m$, and the
  argument for the induction step would finish. Now take
  \[ U_m'=\{f\in U_{m+1}:f\perp f_1\}. \]
  By the same argument, there exists a function $f_2\in U_m'$ which has a zero
  of order $m-1$ at $x$. Using the Taylor expansions
  \[ f_i(r,\omega)=a_ir^{m-1}\cos(m-1)\omega+O(r^m)\qquad i=1,2, \]
  (written in polar coordinates based at $x$, with $a_i\in\c\setminus\{0\}$),
  we see that the
  function $f=a_2f_1-a_1f_2$ is not identically zero, and has a zero of order
  at least $m$ at $x$. This finishes the induction step.
  
  If $f$ has Property~\ref{pro:reallift} then we choose $u=f$. Otherwise, if
  $f$ does not have Property~\ref{pro:reallift} then $\operatorname{Re}[
  \mathcurl Lf]$ is not identically zero, and has a zero of order at least
  $m-1$ at points $y\in\tfold$ such that $\Pi(y)=x$. Using Lemma~\ref{thm:L} we
  see that $u:=\mathcurl L^{-1}(\operatorname{Re}[\mathcurl Lf])$ has
  Property~\ref{pro:reallift}, and has a zero of order at least $m-1$ at $x$.

  \ref{item:unique} For this part we consider the case $m=k+1$ and take any
  point $x\in\Sigma_0$. Part~\ref{item:mdim} shows that there exists a function
  $u_x\in U$ with Property~\ref{pro:reallift} and which has a zero of order at
  least $k$ at $x$. Lemma~\ref{thm:bdzero}
  shows that the zero is of order $k$, and therefore $k$ nodal lines emanate
  from $x$. To prove uniqueness, suppose that $v_x$
  is a linearly independent function which also has a zero of order $k$ at $x$.
  As above, using the Taylor expansions of $u_x$ and $v_x$ at $x$, we may find
  a linear combination of $u_x$ and $v_x$ which is not identically zero and
  which has a zero of order at least $k+1$ at $x$. This contradicts
  Lemma~\ref{thm:bdzero}.

  Due to Lemma~\ref{thm:nodal}, $\Omega\setminus\mathcal N$ is connected, and
  therefore each pair of nodal lines only
  intersect at $x$. The nodal lines must also end at distinct interior boundary
  components.

  Since zeros of order larger than $1$ only occur at points of intersection of
  nodal lines, there can only occur zeros of order $1$ away from $x$. At such
  zeros, the gradient of $u_x$ is non-zero. We may multiply $u_x$ by the local
  gauge $e^{i\phi}$ where $\phi$ is given in equation~\eqref{eqn:Agrad}
  to make it a real valued function. The function $w_x=e^{i\phi}u_x$ has
  locally the same nodal set as $u_x$. Since $w_x$ depends continuously on $x$,
  and the gradient of $w_x$ is non-zero at the nodal set, the nodal lines
  depend smoothly on $x$.

  \ref{item:uni} The proof of this part is similar.\qed
\end{proof}

\begin{proof}[ of Theorem~\ref{thm:main}~\ref{item:mult}]%
  Let $m$ denote the multiplicity of the first eigenvalue of $\op$.
  Lemma~\ref{thm:zero}~\ref{item:mdim} shows that for any point $x\in\partial
  \Omega$ there exists a groundstate of $\op$ which has a zero of order $l\ge m
  -1$ at $x$. Lemma~\ref{thm:bdzero} shows that $l\le k$. This gives the
  universal bound $m\le k+1$, and in particular shows that for $k=1$ we have $m
  \le2$.

  We consider now the case when $k\ge2$ and suppose for a contradiction that
  $m=k+1$. Lemma~\ref{thm:zero}~\ref{item:unique} shows that for each point
  $x$ lying on $\Sigma_0$ there exists a unique eigenfunction $u_x$ which has
  a zero of order $k$ at $x$. Since each $u_x$ has Property~\ref{pro:reallift},
  the nodal set of each $u_x$ \slit{s} $\overline\Omega$. The nodal set of each
  individual $u_x$ has $k$ nodal lines $\{\Gamma_{x,1},\dots,\Gamma_{x,k}\}$,
  emanating from $x$, and each line ends at a distinct interior boundary
  component. We may parametrise each line $\Gamma_{x,i}$ by a path $\gamma_{x,i
  }$ chosen so that $\gamma_{x,i}(0)=x$ and $\gamma_{x,i}(1)\in\Sigma_i$ for
  each $i$. Each path $\gamma_{x,i}$ varies smoothly with $x$.  

  We shall see that if we move $x$ round the boundary $\Sigma_0$, the nodal
  sets of the corresponding functions wind round the holes. After one complete
  turn, we cannot obtain the original nodal set, thus contradicting uniqueness
  of the original eigenfunction. We obtain the contradiction formally as
  follows:
  \doubdiag{%
    \renewcommand\afigwidth{40}%
    \aboundingbox{0}{0}{298}{237}%
    \renewcommand\bfigwidth{40}%
    \bboundingbox{0}{0}{298}{238}}%
    {\fig{movex.eps}{\label{fig:movex}}{%
    \put(76,-8){$x_0$}\put(30,5){$x_t$}\put(105,90){$y_0$}\put(77,97){$y_t$}
    \put(102,47){$\gamma_{x_0,1}$}\put(67,75){$\gamma_{x_t,1}$}
    \put(114,110){$\sigma_1$}\put(195,100){$\sigma_0$}}{}}%
    {\fig{furtherx.eps}{\label{fig:furtherx}}{%
    \put(76,-8){$x_0$}\put(128,-6){$x_t$}\put(88,86){$y_0$}\put(113,99){$y_t$}
    \put(113,116){$\sigma_{1,t}^{-1}$}\put(195,100){$\sigma_{0,t}$}
    \put(100,45){$\gamma_{x_0,1}$}\put(171,67){$\gamma_{x_t,1}^{-1}$}}{}}%

  Let $\sigma_0$ be a closed path which parametrises the outer boundary
  component $\Sigma_0$ of $\Omega$, and which turns once in a clockwise
  direction. For $s\in[0,1]$, let $x_s=\sigma_0(s)$ and let $y_s=\gamma_{x_s,1}
  (1)$. Since $\sigma_0$ is closed, $x_0=x_1$. Also, since $\gamma_{x_s,1}$
  depends
  smoothly on $x_s$, which in turn depends smoothly on $s$, the point $y_s$
  moves smoothly round the inner boundary component $\Sigma_1$. For a fixed
  $t\in[0,1]$ define
  \begin{gather*}
    \sigma_{0,t}(s)=\sigma_0(st)=x_{st}\\
    \sigma_{1,t}(s)=y_{st}
  \end{gather*}
  The paths $\sigma_{0,t}$ and $\sigma_{1,t}$ are parametrisations of segments
  of $\Sigma_0$ and $\Sigma_1$ respectively.
  Note that $\sigma_{0,1}=\sigma_0$ and $\sigma_{1,1}=\sigma_1^p$ for some
  $p\in\z$. For all $t\in[0,1]$ we have
  \begin{equation}\label{eqn:homotopy}
    \sigma_{0,t}^{-1}\circ\gamma_{x_t,1}^{-1}\circ\sigma_{1,t}\circ\gamma_{x_0,
    1}\sim0
  \end{equation}
  where $\circ$ denotes gluing of paths and $\sim$ denotes homotopy. This
  means that the left hand side of~\eqref{eqn:homotopy} is a closed path that
  does not enclose any holes. See Figure~\ref{fig:movex}.
  Setting $t=1$ we get
  \[ \sigma_0^{-1}\circ\gamma_{x_0,1}^{-1}\circ\sigma_1^p\circ\gamma_{x_0,1}
     \sim0, \]
  and therefore
  \[ \sigma_1^p\sim\gamma_{x_0,1}^{-1}\circ\sigma_1^p\circ\gamma_{x_0,1}\sim
     \sigma_0. \]
  This gives us a contradiction because the path $\sigma_1^p$ is not homotopic
  to $\sigma_0$. Hence $m\le k$.

  Finally we consider the case where $k$ is even and $k\ge4$. Let $\overline
  \Omega'$ denote the closure $\overline\Omega$ of our region with the points
  of the outer boundary identified.
  Let $D_{k-1}\subset\r^2$ denote an open disk with $k-1$ smaller, disjoint,
  closed disks removed. There exists a homeomorphism
  \begin{equation}\label{eqn:X}
    X:\overline\Omega'\rightarrow{\overline D}_{k-1}
  \end{equation}
  \singdiag{%
  \renewcommand\afigwidth{90}%
  \aboundingbox{0}{0}{499}{183}}{%
  \fig{evenodd.eps}{\label{fig:evenodd}}{%
  \put(171,60){$p$}%
  \put(102,29){$X$}%
  \put(43,8){$\overline\Omega$}%
  \put(165,8){${\overline D}_{k-1}$}%
  \put(79,59){$\Sigma_0$}%
  \put(68,25){$\Sigma_1$}%
  }{}}%
  such that $X$ restricted to
  $\Omega$ is smooth, and such that the boundary component $\Sigma_1$ maps to
  the outer boundary of ${\overline D}_{k-1}$. See Figure~\ref{fig:evenodd}.
  One can imagine $X$ as a composition of mapping $\overline\Omega$ onto the
  surface of a sphere, deforming it so that $\Sigma_1$ becomes very large and
  $\Sigma_0$ very small, and then finally pulling off the sphere. Let $p:=X(
  \Sigma_0)\in D_{k-1}$, so that $X(\Omega)=D_{k-1}\setminus\{p\}$.

  Let $\mathcal N$ be a set which \slit{s} $\overline\Omega$. We claim that
  $X(\mathcal N)$ \slit{s} ${\overline D}_{k-1}$. For since $k$ is even, the
  number
  of nodal lines hitting the outer boundary component $\Sigma_0$ is even
  (possibly zero). This corresponds to an even number of paths in $X(\mathcal N
  )$ starting or finishing at $p$. These paths can be paired together to link
  distinct boundary components. Since $X^{-1}$ is a smooth bijection
  away from $p$, the resulting paths are still piecewise smooth. It is easy to
  verify that all the other \slit{ting} conditions are satisfied.

  Suppose for a contradiction that $m=k$. For $s\in[0,1]$, let $x_s=\sigma_1(s
  )$ be a point on the interior boundary component $\Sigma_1$ of $\Omega$.
  Lemma~\ref{thm:zero}~\ref{item:uni} shows that there exists a unique $u_{x_s}
  \in U$ (up to multiplication by a complex constant) which has a zero of order
  $k-1$ at $x_s$. The nodal set $\mathcal N(u_{x_s})$ consists of $k-1$ nodal
  lines emanating from $x_s$. As shown above, the set $S_s:=X(\nod[u_{x_s}])$
  \slit{s} ${\overline D}_{k-1}$ and consists of $k-1$ lines emanating from the
  point $y_s=X(x_s)$ on the outer boundary of ${\overline D}_{k-1}$.

  We have thus constructed a family of \slit{ting} sets $S_s$ which depends
  continuously on the parameter $s\in[0,1]$, and such that $S_0=S_1$. By moving
  the point $y_s$ round the outer boundary of ${\overline D}_{k-1}$ and using
  the homotopy argument above, we obtain a similar contradiction. Hence
  $m\le k-1$.\qed
\end{proof}

\begin{proof}[ of Theorem~\ref{thm:main}~\ref{item:emptynodal}]
  Suppose that $k=1$ and that the multiplicity of the first eigenvalue is two.
  Suppose for a contradiction that there exist two linearly independent
  groundstates $v_1$ and $v_2$ such that the set $S=\nod[v_1]\cap\nod[v_2]$ is
  non-empty, and let $z$ be any point in $S$. Since $\{v_1,v_2\}$ is a basis
  of the groundstate eigenspace $U$ of $\op$, the nodal set of
  every function $u\in U$ contains the point $z$.

  From Lemma~\ref{thm:zero}~\ref{item:unique} we see that for each point $x$ on
  the outer boundary $\Sigma_0$ of $\Omega$ there exists a unique eigenfunction
  $u_x\in U$
  such that $x\in\nod[u_x]$. If we start $x$ at the point $x_0=\sigma_0(0)$ and
  then move $x$ continuously round the outer boundary $\Sigma_0$ once in an
  clockwise direction then the segment of the nodal line joining $x$ to $z$
  deforms continuously and winds around the inner boundary $\Sigma_1$ (see
  Figure~\ref{fig:nodal}). The resulting nodal line is
  \singdiag{%
    \renewcommand\afigwidth{40}%
    \bboundingbox{0}{0}{282}{224}}{%
    \fig{contrad.eps}{\label{fig:nodal}}{%
    \put(25,150){$x_0$}%
    \put(49,112){$z$}}{}}%
  different to the original, thus contradicting uniqueness of the eigenfunction
  $u_{x_0}$. This argument can be formalised using a homotopy argument similar
  to that found in the proof of part~\ref{item:mult}.

  Suppose that $u=\alpha_1u_1+\alpha_2u_2$ where $\alpha_1\overline{\alpha_2}
  \not\in\r$. Since each function $\mathcurl Lu_i$ is real valued (see the
  construction of the $u_i$ in the proof of
  Theorem~\ref{thm:main}~\ref{item:cutting}), we have
  \[ \nod[\mathcurl L(\overline{\alpha_2}u)]=\nod[\alpha_1\overline{\alpha_2}
     \mathcurl Lu_1+|\alpha_2|^2\mathcurl Lu_2]=\nod[\mathcurl Lu_1]\cap\nod[
     \mathcurl Lu_2]. \]
  Since the nodal sets of $u_1$ and $u_2$ do not intersect, we have
  \begin{equation*}
    \nod[u]=\Pi(\nod[\mathcurl Lu])\subseteq\Pi(\nod[\mathcurl Lu_1])\cap\Pi(
    \nod[\mathcurl Lu_2])=\nod[u_1]\cap\nod[u_2]=\emptyset.
  \end{equation*}
  For the case $k=2$, the proof uses the map $X:\overline\Omega'\rightarrow
  \overline D_1$ (see equation~\eqref{eqn:X}) to essentially reduce the region
  with two holes to the single hole case.\qed
\end{proof}

\begin{proof}[ of inequality~\eqref{eqn:fourth} from
               Theorem~\ref{thm:properties}]%
  Suppose that $k=1$, and let $A_1$ and $A_2$ be magnetic vector potentials,
  where $A_1$ has circulation $1/2$. Let $\Phi$ denote the circulation of
  $A_2$. Suppose for a contradiction
  that $\Phi\not\in1/2+\z$ and that $\lambda_1(\op[A_2])\ge\lambda_1(\op[A_1])
  $. Using Theorem~\ref{thm:main}~\ref{item:cutting}, there exists a
  groundstate $u_1$ of $\op[A_1]$ which has a nodal set $\mathcal N$
  which \slit{s} $\overline\Omega$. As we are in the single hole case, the
  nodal set consists of a single line $\Gamma$ which joins the outer boundary
  to the inner boundary.

  We shall need an operator $\op[\Gamma,A_2]$, which has extra Dirichlet
  boundary conditions imposed along the line $\Gamma$. This is defined formally
  as the self-adjoint operator corresponding to the restriction of the closed
  quadratic form $\qf[A_2]$ (defined in~\eqref{eqn:quaddef}) to the domain
  \begin{equation*}
    \q_\Gamma^{\operatorname{Neu}}=\{u\in\q^{\operatorname{Neu}}=W^{1,2}(\Omega
    ):u|_{\Gamma}=0\}.
  \end{equation*}
  Using our supposition, and the fact that the nodal set of $u_1$ consists of
  the line $\Gamma$, we have
  \begin{equation}\label{eqn:extraD1}
    \lambda_1(\op[A_2])\ge\lambda_1(\op[A_1])=\lambda_1(\op[\Gamma,A_1])
  \end{equation}
  Since $\Omega\setminus\Gamma$ is simply connected, $\op[\Gamma,A_1]$
  is unitarily equivalent to $\op[\Gamma,A_2]$, and therefore
  \begin{equation}\label{eqn:extraD2}
    \lambda_1(\op[\Gamma,A_1])=\lambda_1(\op[\Gamma,A_2])
    =\inf_{u\in{\q^{\operatorname{Neu}}_\Gamma}(\Omega)}\qf[A_2](u)\ge\lambda_1
    (\op[A_2]).
  \end{equation}
  We have equality in~\eqref{eqn:extraD1} and~\eqref{eqn:extraD2}, and
  therefore the groundstate $u_2\in\q^{\operatorname{Neu}}$ of
  $\op[\Gamma,A_2]$ is also a groundstate of $\op[A_2]$
  The nodal sets of $u_1$ and $u_2$ both contain $\Gamma$.

  Since $\curl A_1=\curl A_2=0$ in the connected set $\Omega\setminus\Gamma$,
  there exist smooth functions $\phi_1,\phi_2:\Omega\setminus\Gamma
  \rightarrow\r$ such that
  $\grad\phi_i=A_i$. The functions $\phi_1$ and $\phi_2$ supply us with gauge
  transformations $e^{i\phi_1}$ and $e^{i\phi_2}$, from which we see both
  $e^{i\phi_1}u_1$ and $e^{i\phi_2}u_2$ are groundstates of $\op[\Gamma,0]
  $. By uniqueness of the groundstate of a non magnetic \schroedinger operator,
  we have
  \[ u_2=\lambda e^{i(\phi_2-\phi_1)}u_1 \]
  for some constant $\lambda\in\c\setminus\{0\}$. Let $\phi_3=\phi_2-\phi_1$.
  Since both $u_1$ and $u_2$ are smooth functions on $\Omega$ we may extend
  $\phi_3$ to a $C^1$ multivalued function on $\Omega$. The values
  that $\phi_3$ takes at a point differ by multiples of $2\pi$.
  Hence for a path $\sigma$ which circulates $\Omega$ once
  \[ \frac1{2\pi}\int_\sigma A_2\cdot\d\x=\frac1{2\pi}\int_\sigma A_1\cdot\d\x+
     \frac1{2\pi}\int_\sigma(A_2-A_1)\cdot\d\x=\frac12+\frac1{2\pi}\int_\sigma\d
     \phi_3=\frac12+l. \]
  This contradicts our assumption that $\Phi\not\in1/2+\z$.\qed
\end{proof}

\begin{remarks}\begin{enumerate}
  \item Using semiclassical arguments as in~\cite{Helf1}, we can show that for
  $k\ge2$, the first eigenvalue is not necessarily maximised for circulation
  $(1/2,\dots,1/2)$. However, we may use methods similar to those in the above
  proof to show that
  \begin{equation}\label{eqn:inf}
    \lambda_1(1/2,\dots,1/2)=\inf_{S\in\mathcurl S}\lambda_1(\op[S,0]),
  \end{equation}
  where $\mathcurl S$ is the collection of all sets $S$ which \slit{} $\Omega$,
  and
  where $\op[S,0]$ is defined (as in the above proof) to have extra Dirichlet
  boundary conditions along $S\in\mathcurl S$.
  \item

  In~\cite{BergRubi}, Berger and Rubinstein investigate the single hole case,
  in which a nodal set $S$ which \slit{s} $\overline\Omega$ is simply a curve
  $\Gamma$ joining the outer boundary to the inner boundary. As a condition for
  their analysis, they require that $\Omega$ has the property that the
  minimum~\eqref{eqn:inf} is attained for a finite number of curves $S$.
  
  Theorem~\ref{thm:main} shows that this requirement is in fact misleading: If
  the groundstate of $\op$ is simple then the infimum~\eqref{eqn:inf} is
  attained for one curve $S\in\mathcurl S$. Otherwise, the multiplicity of
  the first eigenvalue is two, and the infimum is attained for infinitely many
  curves.
  \end{enumerate}
\end{remarks}

\section{Additional results and examples}

\begin{proposition}\label{thm:boundlines}
  If a collection of paths $\{\gamma_1,\dots,\gamma_n\}$ \slit{s} a region
  $\Omega$ with $k$ holes then $k/2\le n\le k$.
\end{proposition}

\begin{proof}
  The lower bound on $n$ is elementary because there are an odd number of lines
  (i.e. at least one) leaving each of the $k$ holes. There must therefore
  be at least $k/2$ lines.

  We finally prove the upper bound on $n$.
  Let $\sigma_0$ be a closed path which parametrises the outer boundary
  $\Sigma_0$ of $\Omega$, and let $\sigma_1,\dots,\sigma_k$ be closed paths
  which parametrise the $k$ other boundary components $\Sigma_1,\dots,\Sigma_k
  $. Define
  \begin{gather}
  S_0=\bigcup_{i=0}^k\sigma_i(0)\label{eqn:0cells}\\
  S_1=\left(\bigcup_{i=0}^k\big\{\sigma_i\big((0,1)\big)\big\}\right)\cup\left(
  \bigcup_{j=1}^n\big\{\gamma_j\big([0,1]\big)\big\}\right)\label{eqn:1cells}\\
  S_2=\{\Omega\setminus\mathcal N\}.
  \end{gather}
  
  Let $(N_0,N_1,N_2)=(k+1,k+1+n,1)$ be the triple of integers associated to
  this decomposition, in which $N_i$ is the number of elements in the
  collection $S_i$. The decomposition $D$ is not a standard CW decomposition of
  $\overline\Omega$, and
  therefore the number $N:=N_0-N_1+N_2$ will not yield the Euler number
  $\chi(\Omega)=-k+1$. It is however possible to modify the decomposition
  to make it into a proper CW decomposition in two steps:
  \begin{enumerate}
    \item We first add vertices where intersections of elements of $S_1$
    occur at points which are not in $S_0$.
    This step will decompose some elements of $S_1$ into smaller parts but
    leaves the element $\Omega\setminus\mathcal N$ of $S_2$ unaltered. Let
    $S_0'$ denote the new collection of vertices.
    \item If $\Omega\setminus\mathcal N$ is not simply
    connected then the second step is to add some extra lines, which
    begin and end at already existing vertices in $S_0'$ in order to break
    up (without disconnecting) the region into a single simply connected
    $2$-cell.
  \end{enumerate}
  Note that after each step, $S_2'$ still consists of just one connected open
  set, so $N_2'=1$, whilst the number $N_0'-N_1'$ of vertices minus lines
  does not increase. It follows that
  \[ N_0-N_1+N_2\ge N_0'-N_1'+N_2'=\chi(\Omega). \]
  Substituting in $N_0=k+1$, $N_1=k+1+n$, $N_2=1$, and $\chi(\Omega)=-k+1$, we
  obtain
  \begin{equation*}
    n\le k.\tag*\qed
  \end{equation*}
\end{proof}

\begin{example}
  The example of the circle $S^1$ is interesting to analyse. Consider the
  operator
  \[ P_\alpha=-(\partial_\phi-i\alpha)^2. \]
  on $L^2(S^1)$. The spectrum can be easily seen to be
  \[ \sigma(P_\alpha)=\{(n-\alpha)^2:n\in\z\}, \]
  and therefore
  \[ \lambda_1(P_\alpha)=\min_{n\in\z}(n-\alpha)^2. \]
 
  When $\alpha$ is an integer, the first eigenvalue is $0$ and is simple and
  the corresponding eigenfunction is $\exp i\alpha\phi$. The first eigenvalue
  is actually simple whenever $\alpha$ is not a half-integer.

  On the other hand, if $\alpha$ is a half-integer, the first eigenvalue is
  $1/4$, with multiplicity two. The corresponding eigenspace is spanned by the
  functions $1$ and $\exp(i\phi)$ (or alternatively by the functions $\exp(i
  \phi/2)\cos(\phi/2)$ and $\exp(i\phi/2)\sin(\phi/2$)), and one can
  parametrise all the resulting eigenfunctions, in terms of a parameter $\phi_
  0$, by $\exp(i\phi/2)\sin((\phi-\phi_0)/2)$.

  It is easy to see how the degeneracy of the first eigenvalue disappears when
  considering
  \[ P_{\alpha,\epsilon,v}=-(\partial_\phi-i\alpha)^2+\epsilon v(\phi), \]
  perturbatively as $\epsilon\neq0$ is small, provided $v(\phi)$ satisfies the
  condition
  \[ \int_0^{2\pi}v(\phi)e^{i\phi}\neq0. \]
\end{example}

\begin{example}\label{exa:sharp}
In \cite[Subsection 7.3]{Helf2}, an example is given in which the multiplicity
of the first eigenvalue is two. The domain $\Omega$ and potential $V$ are
symmetric under the map $S:z\mapsto-z$, and the magnetic potential is given
explicitly by
\[ A=\frac\Phi{2\pi r^2}\binom{-y}x. \]
If we take the case when the flux is an half-integer and 
we compose the operator $K$ (see Remark~\ref{rem:K}) with the operator $S$
defined by
\[ (Su)(z)=u(Sz), \]
the operator 
\[ M=SK \]
commutes with $P_{A,V}$ and satisfies
\[ M^2=-I. \]
Kramer's theorem shows that the multiplicity is at least two. One
can indeed show that $u$ and $Mu$ are linearly independent.

An alternative proof is simply to say that $Su$ is also an
eigenvector with nodal set $S\gamma$, where $\gamma$ is the nodal set of $u$.
\end{example}

\subsection*{Acknowledgments}

The authors wish to thank L.~Friedlander, P.~Michor for many interesting
discussions and J.~Rubinstein for useful correspondence. B.~Helffer and
T.~Hoffmann-Ostenhof are grateful to A.~Laptev for inviting them to Stockholm,
where this research was initiated.

\bibliographystyle{amsalpha}
%\ignore{
\providecommand{\bysame}{\leavevmode\hbox to3em{\hrulefill}\thinspace}

%} % end of \ignore{...
%\bibliography{/.automount/Keen/users/mowen/Tex/Bibliog/mybib}

\begin{thebibliography}{HOHON}

\bibitem[BCC98]{BessColbCour}
G.~Besson, B.~Colbois, and G.~Courtois, \emph{Sur la {multiplicit\'e} de la
  {premi\`ere} valeur propre de {l'op\'erateur} de {Schr\"odinger} avec champ
  {magn\'etique} sur la {sph\`ere} {$S^2$}}, Trans. Amer. Math. Soc.
  \textbf{350} (1998), 331--345.

\bibitem[Ber55]{Bers}
L.~Bers, \emph{Local behaviour of solutions of general linear equations},
  Commun. Pure Appl. Math. \textbf{8} (1955), 473--496.

\bibitem[BR97]{BergRubi}
J.~Berger and J.~Rubinstein, \emph{On the zero set of the wave function in
  superconductivity}, Preprint {mp\_arc} 97, 1997.

\bibitem[Che76]{Chen}
S.~Y. Cheng, \emph{Eigenfunctions and nodal sets}, Commentarii. Math. Helv.
  \textbf{51} (1976), 43--55.

\bibitem[{Col}93]{Coli}
Y.~{Colin de Verdi\`ere}, \emph{Multiplicit\'es des valeurs propres.
  {Laplaciens} discrets et laplaciens continus}, Rend. Mat. Appl. \textbf{VII}
  (1993), 433--460.

\bibitem[EMQ94]{ElliMataQi}
C.~M. Elliott, H.~Matano, and T.~Qi, \emph{Zeros of complex {G}inzburg-{L}andau
  order parameter with applications to superconductivity}, Eur. J. Appl. Math.
  \textbf{5} (1994), 431--448.

\bibitem[GHL90]{GallHuliLafo}
S.~Gallot, D.~Hulin, and J.~Lafontaine, \emph{Riemannian geometry}, 2nd ed.,
  Universitext, Springer-Verlag, 1990.

\bibitem[GT83]{GilbTrud}
N.~Gilbarg and S.~Trudinger, \emph{Elliptic partial differential equations of
  second order}, 2nd ed., Grundlehren der mathematischen {W}issenschaften, no.
  224, Springer-Verlag, 1983.

\bibitem[Hel88a]{Helf1}
B.~Helffer, \emph{Effet {d'Aharonov-Bohm} sur un \'etat born\'e de l'\'equation
  de {Schr\"odinger}}, Commun. Math. Phys. \textbf{119} (1988), 315--329.

\bibitem[Hel88b]{Helf2}
B.~Helffer, \emph{Semi-classical analysis for the {Schr\"odinger} operator and
  applications}, Lecture notes in mathematics, no. 1336, Springer-Verlag, 1988.

\bibitem[Hel94]{Helf3}
B.~Helffer, \emph{On spectral theory for {Schr\"odinger} operators with
  magnetic potentials}, Adv. Stud. Pure Math. \textbf{23} (1994), 113--141.

\bibitem[HOHON]{HoffHoffNadi}
M.~Hoffmann-Ostenhof, T.~Hoffmann-Ostenhof, and N.~Nadirashvili, \emph{On the
  multiplicity of eigenvalues of the {Laplacian} on surfaces}, Submitted.

\bibitem[HOMN]{HoffMichNadi}
T.~Hoffmann-Ostenhof, P.~Michor, and N.~Nadirashvili, \emph{Bounds on the
  multiplicity of eigenvalues for fixed membranes}, Preprint 1998.

\bibitem[Kos80]{Kosn}
C.~Kosniowski, \emph{A first course in algebraic topology}, Cambridge
  University Press, 1980.

\bibitem[LO77]{LaviOCar}
R.~Lavine and M.~O'Carroll, \emph{Ground state properties and lower bounds for
  energy levels of a particle in a uniform magnetic field and external
  potential}, J. Math. Phys. \textbf{18} (1977), 1908--1912.

\bibitem[LP62]{LittPark}
W.A. Little and R.D. Parks, \emph{Observation of quantum periodicity in the
  transition temperature of a superconducting cylinder}, Phys. Rev. Lett.
  \textbf{9} (1962), 9--12.

\bibitem[Nad88]{Nadi}
N.~S. Nadirashvili, \emph{Multiple eigenvalues of the {Laplace} operator},
  Math. USSR, Sb. \textbf{61} (1988), 225--238.

\bibitem[Sim79]{Simo}
B.~Simon, \emph{Kato's inequality and the comparison of semigroups}, J. Funct.
  Anal. \textbf{32} (1979), 97--101.

\bibitem[Wlo82]{Wlok}
J.~Wloka, \hspace{1em}\emph{Partielle {D}ifferentialgleichungen. {S}obolevr{\"a}ume und
  {R}and\-wert\-auf\-gaben}, B. G. Teubner, 1982.

\end{thebibliography}

\small
\noindent{\scshape
Bernard Helffer, D\'epartement de Math\'ematiques, B\^atiment 425, Universit\'e
Paris-Sud, F-91405 Orsay C\'edex, France.}\\
{\itshape E-mail address: }{\sffamily bernard.helffer@@math.u-psud.fr}\\

\noindent{\scshape
Maria Hoffmann-Ostenhof, Institut f\"ur Mathematik, Universit\"at\\ Wien,
Strudlhofgasse 4, A-1090 Wien, Austria.}\\
{\itshape E-mail address: }{\sffamily mho@@nelly.mat.univie.ac.at}\\

\noindent{\scshape
Thomas Hoffmann-Ostenhof, Insitute f\"ur Theoretische Chemie, Universit\"at
Wien, W\"ahringerstra\ss{}e 17, A-1090 Wien, Austria.}\\
{\itshape E-mail address: }{\sffamily hoho@@itc.univie.ac.at}\newline
\textrm{or}\\
{\scshape Thomas Hoffmann-Ostenhof, Erwin \schroedinger International Insitute
for Mathematical Physics (ESI), Botzmanngasse 9, A-1090 Wien, Austria.}\\
{\itshape E-mail address: }{\sffamily thoffman@@esi.ac.at}\newline

\noindent{\scshape
Mark Owen, Erwin \schroedinger International Insitute for
Mathematical Physics (ESI), Botzmanngasse 9, A-1090 Wien, Austria.}\\
{\itshape E-mail address: }{\sffamily mowen@@esi.ac.at}\newline

\end{document}